\documentclass{article}

\usepackage{amsfonts}
\usepackage{psfrag}
\usepackage{graphics}
\usepackage{amsmath}
\usepackage{amsthm}
\usepackage{amssymb}
\usepackage{epsfig}

\newtheorem{theorem}{Theorem}
\newtheorem{corollary}[theorem]{Corollary}
\newtheorem{lemma}[theorem]{Lemma}

\theoremstyle{definition}
\newtheorem{definition}[theorem]{Definition}
\theoremstyle{remark}

\theoremstyle{definition}

\newfont{\valami}{ptmr8r scaled 1200}
\newfont{\kisvalami}{ptmr8r scaled 1000}

\newcommand{\rank}{\operatorname{rank}}

\begin{document}

\author{
Leonidas Pitsoulis  \\
Aristotle University of Thessaloniki
\and
Konstantinos Papalamprou \\
London School of Economics
\and
Gautam Appa\\
London School of Economics
\and
Bal\'azs Kotnyek \\
I3S {\&} Inria Sophia Antipolis
}
\title{On the representability of totally unimodular matrices on bidirected graphs}

\maketitle

\begin{abstract}
Seymour's famous decomposition theorem for regular matroids states that any totally unimodular (TU) matrix can be 
constructed through a series of composition operations called $k$-sums starting from network matrices and their 
transposes and two  compact representation matrices $B_{1}, B_{2}$ of a certain ten element matroid. Given that $B_{1}, B_{2}$ 
are binet matrices we examine the $k$-sums of network and binet matrices. It is shown that the $k$-sum of a network and 
a binet matrix is a binet matrix, but binet matrices are not closed under this operation for $k=2,3$. A new class of 
matrices is introduced the so called {\em tour matrices}, which generalises network,  binet and totally unimodular matrices. 
For any such matrix there exists a bidirected graph such that the columns represent a collection of closed tours in the graph. 
It is shown that tour matrices are closed under $k$-sums, as well as under pivoting and other elementary operations on its rows 
and columns.  
Given the constructive proofs of the above results regarding the $k$-sum operation and existing recognition algorithms for network and
binet matrices, an algorithm is presented which constructs a bidirected graph for any TU matrix. 
\end{abstract}


% \setcounter{tocdepth}{4}
% \tableofcontents

%%%%%%%%%%%%%%%%%%%%%%%%%%%%%%%%%%%%%%%%%%%%%%%%%%%%%%%%%%%%%%%%%%%%%%%%%%%%%%%%%%%%%%%%%%

\section{Introduction}
Totally unimodular matrices are a class of $\{0,\pm1\}$ matrices which is of great importance to 
combinatorial optimisation since they describe a special class of polynomial
time solvable integer programs. Specifically, every integer program which is defined by a totally unimodular constraint 
matrix can be solved as a linear program by relaxing the integrality constraint since the associated polyhedron 
is integral. Although there exist various equivalent characterisations for this class of matrices, it was Seymour's 
decomposition theory~\cite{Seymour:1980} developed for the associated regular matroids, that yielded a polynomial time 
algorithm for recognising them. 
Seymour's decomposition theorem states that all totally unimodular matrices can be constructed 
recursively by applying $k$-sum operations $(k=1,2,3)$ on network matrices, their transposes and two 
totally unimodular matrices $B_1$ and $B_2$. These sum operations, are essentially matrix 
operations which preserve certain structural properties. 
Combined with the fact that the matrices $B_1$ and $B_2$ are easily recognisable, and Tutte's theory 
for recognising network matrices, Seymour's theorem implies an algorithm to check whether a given 
matrix is totally unimodular or not. Moreover, and maybe even more importantly, it also provides a 
framework for graphical representation of totally
unimodular matrices. 
% It is partial in the sense that $B_1$ and $B_2$ as well as the transposes
% of network matrices which correspond to non planar graphs, are not network matrices therefore
% do not have an associated directed graph. 
Bidirected graphs are a generalisation of directed graphs, and can be represented algebraically by the 
so-called binet matrices in the same way network matrices represent directed graphs. 
Appa and Kotnyek~\cite{AppKot:2005} have shown that $B_1$ and $B_2$ can be represented on 
bidirected graphs since they have been proved to be binet. Since bidirected graphs generalise directed graphs, all the building 
blocks of totally unimodular matrices or their transposes are representable on bidirected graphs. 
In this work we show that every totally unimodular matrix has an associated bidirected graph representation, 
which provides a partial interpretation of the nice integrality property of the associated polyhedron and 
may provide the means of devising a combinatorial algorithm for solving the related integer programming problem. 

Initially we show constructively that the $k$-sum of two network matrices is a network matrix and 
that of a network and a binet matrix is a binet matrix. 
% Moreover, given the network or binet 
% representations of network or binet matrices respectively, we provide graphical methods to construct 
% representations of the associated operations on these matrices. 
However, for $k=2,3$ we show that 
the $k$-sum of two binet matrices is not necessarily a binet matrix. 
Based on this we can state that not all totally unimodular matrices are binet. 
To pursue graphical representability further a new class of $\{0,\pm 1\}$ matrices is introduced, the so-called {\em tour matrices}, 
which represent closed tours on bidirected graphs. We show that network matrices as well as 
$B_1$ and $B_2$ are tour matrices, and in contrast to binet matrices, it is also shown that
tour matrices are closed under $k$-sums. 
This means that totally unimodular 
matrices not previously associated with bidirected graphs can now be represented on bidirected graphs.  

The paper is organised as follows. Section~\ref{sec_prel} presents all the preliminary theory
regarding network matrices, bidirected graphs and binet matrices, totally unimodular matrices as
well as the definition of the $k$-sum operations. 
In section~\ref{sec_ksums} we examine the operation of $k$-sums of network and binet matrices, 
where the most general case for $k=3$ is treated and a graphical construction of the operation
is presented. The negative result in this section is that binet matrices are not closed under
$k$-sums. Tour matrices are defined in section~\ref{subsec_tour_properties} where various properties are
proved. In section~\ref{subsubsec_ksums_tour} we show that tour matrices are closed under $k$-sums, while
in section~\ref{subsubsec_algorithm} we gather all the results presented on the paper on an algorithm
for constructing a bidirected graph of any TU matrix.

%%%%%%%%%%%%%%%%%%%%%%%%%%%%%%%%%%%%%%%%%%%%%%%%%%%%%%%%%%%%%%%%%%%%%%%%%%%%%%%%%%%%%%%%%%

\section{Preliminaries}\label{sec_prel}

%%%%%%%%%%%%%%%%%%%%%%%%%%%%%%%%%%%%%%%%%%%%%%%%%%%%%%%%%%%%%%%%%%%%%%%%%%%%%%%%%%%%%%%%%%
\subsection{Graphs and Network Matrices} \label{sec_dirgraphs}

A \emph{directed graph} $G(V,E)$ consists of a finite set of nodes $V$ and a family $E$ of 
\emph{ordered} pairs of $V$. For an edge $e=(u,v)$, $u$ and $v$ are called the \emph{end-nodes} of 
$e$; $u$ is called the \emph{tail} of $e$ and $v$ the \emph{head} of $e$. We also say that 
$e=(u,v)$ \emph{leaves} $u$ and \emph{enters} $v$.
The \emph{node-edge incidence matrix}  of a directed graph $G(V,E)$ is the $V\times{E}$ matrix 
$D_{G}$ with
\[
D_{G}(v,e)=
\left\{ 
\begin{array}{rl}
-1 & \textrm{if $v$ is a tail of $e$} \\
+1  & \textrm{if $v$ is a head of $e$} \\
0  & \textrm{otherwise,}
\end{array} \right.
\]\\
for any $v\in{V}$ and any non-loop $e\in{E}$. If $e$ is a loop, we set $D_{G}(v,e):=0$ for 
each vertex $v$.
The definition for the network matrices goes as follows:
\begin{definition}
Let $D_{G}=[R|S]$ be the incidence matrix of a directed graph $G(V,E)$ minus an arbitrary row, where $R$ is a basis of the 
column space of $D_{G}$. The matrix $N_{G}=R^{-1}S$ is called a \emph{network matrix}.
\end{definition}
For material related to graphs and network matrices the reader is referred to~\cite{Schrijver:2004}. 
% The following are well known properties of network matrices, which can be shown easily since they
% do not affect linear independence:
% \begin{itemize}
% \item [(a)] Network matrices are closed under row and column deletions or duplications.
% \item [(b)] Network matrices are closed under addition of a unitary (i.e. with one element equal to $1$ and all other $0$) row or column.
% \item [(c)] Network matrices are closed under multiplication of a row or column by $-1$.
% \item [(d)] Network matrices are closed under pivoting (in $\mathbb{R}$).
% \item [(e)] If $N$ is a network matrix then $[I~N]$ is a network matrix also, where $I$ is the 
% identity matrix.  
% \end{itemize}

%%%%%%%%%%%%%%%%%%%%%%%%%%%%%%%%%%%%%%%%%%%%%%%%%%%%%%%%%%%%%%%%%%%%%%%%%%%%%%%%%%%%%%%%%%
\subsection{Bidirected Graphs and Binet Matrices}
A bidirected graph $\Sigma(V,E)$ is defined over a finite node set $V$ and an edge set 
$E\subseteq V\times V$. There are four types of edges: a \emph{link} has two different end-nodes, 
a \emph{loop} has two end-nodes  that coincide, a \emph{half-edge} has one end-node, 
and a \emph{loose edge} which has no end-nodes \cite{AppKot:2005}.

\begin{figure}[ht]
\begin{center}
\centering
\psfrag{r1}{\footnotesize $r_1$}
\psfrag{r2}{\footnotesize $r_2$}
\psfrag{r3}{\footnotesize $r_3$}
\psfrag{r4}{\footnotesize $r_4$}
\psfrag{r5}{\footnotesize $r_5$}
\psfrag{r6}{\footnotesize $r_6$}
\psfrag{r7}{\footnotesize $r_7$}
\psfrag{r8}{\footnotesize $r_8$}
\psfrag{s1}{\footnotesize $s_1$}
\psfrag{s2}{\footnotesize $s_2$}
\psfrag{s3}{\footnotesize $s_3$}
\psfrag{s4}{\footnotesize $s_4$}
\psfrag{s5}{\footnotesize $s_5$}
\psfrag{s6}{\footnotesize $s_6$}
\psfrag{v1}{\footnotesize $v_1$}
\psfrag{v2}{\footnotesize $v_2$}
\psfrag{v3}{\footnotesize $v_3$}
\psfrag{v4}{\footnotesize $v_4$}
\psfrag{v5}{\footnotesize $v_5$}
\psfrag{v6}{\footnotesize $v_6$}
\psfrag{v7}{\footnotesize $v_7$}
\psfrag{v8}{\footnotesize $v_8$}
\hfill
\includegraphics*[scale=0.6]{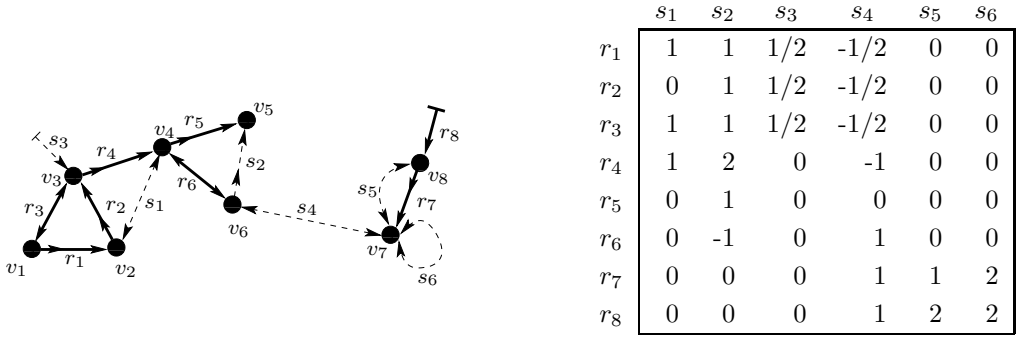}\hfill
\raisebox{10ex}{\ \ \ \ \begin{tabular}{r|rrrrrr|} 
\multicolumn{1}{c}{} & \multicolumn{1}{c}{$s_1$} &\multicolumn{1}{c}{$s_2$} & 
\multicolumn{1}{c}{$s_3$} & \multicolumn{1}{c}{$s_4$} & \multicolumn{1}{c}{$s_5$} & \multicolumn{1}{c}{$s_6$}\\ \cline{2-7} 
$r_1$ & 1 & 1 & 1/2 & -1/2& 0 & 0\\  
$r_2$ & 0 & 1 & 1/2 & -1/2& 0 & 0\\  
$r_3$ & 1 & 1 & 1/2 & -1/2& 0 & 0\\  
$r_4$ & 1 & 2 & 0   & -1  & 0 & 0\\  
$r_5$ & 0 & 1 & 0   &  0  & 0 & 0\\  
$r_6$ & 0 &-1 & 0   &  1  & 0 & 0\\  
$r_7$ & 0 & 0 & 0   &  1  & 1 & 2\\  
$r_8$ & 0 & 0 & 0   &  1  & 2 & 2\\ \cline{2-7} 
\end{tabular}}\hspace*{\fill}
\end{center}
\caption{An example of a binet graph, and its binet matrix. } \label{fig:fundex}
\end{figure}

Every edge is assigned a \emph{sign}, so that half-edges are always negative; loose edges are always 
positive; links and loops can be positive or negative. The edges are \emph{oriented}, i.e., we label 
the end-nodes of the edges by $+1$ or $-1$. The labels of a positive edge are different, those of a 
negative edge are the same. If an end-node of an edge is labeled with $+1$, then it is an \emph{in-node} 
of the edge, otherwise an \emph{out-node}. These names come from the  graphical representation of 
bidirected graphs, where incoming and outgoing arrows on an edge represent positive and  negative 
labels. For example in the bidirected graph shown in Figure~\ref{fig:fundex}, 
edge $r_1$ is a positive link;  $r_3$ is a negative; $s_6$ is a negative loop; and $r_8$ is a half-edge. 
Loose edges and positive loops are not depicted in this illustration.  
%All negative edges will be called \emph{bidirected} while all positive edges \emph{directed}. 
%\begin{figure}
%\begin{center}
%\centering
%\psfrag{r1}{\footnotesize $e_1$}
%\psfrag{r2}{\footnotesize $e_2$}
%\psfrag{r3}{\footnotesize $e_3$}
%\psfrag{s3}{\footnotesize $e_4$}
%\psfrag{r4}{\footnotesize $e_5$}
%\psfrag{s1}{\footnotesize $e_6$}
%\psfrag{r5}{\footnotesize $e_7$}
%\psfrag{r6}{\footnotesize $e_8$}
%\psfrag{s2}{\footnotesize $e_9$}
%\psfrag{s4}{\footnotesize $e_{10}$}
%\psfrag{s5}{\footnotesize $e_{11}$}
%\psfrag{s6}{\footnotesize $e_{12}$}
%\psfrag{r7}{\footnotesize $e_{13}$}
%\psfrag{r8}{\footnotesize $e_{14}$}
%\psfrag{v1}{\footnotesize $v_1$}
%\psfrag{v2}{\footnotesize $v_2$}
%\psfrag{v3}{\footnotesize $v_3$}
%\psfrag{v4}{\footnotesize $v_4$}
%\psfrag{v5}{\footnotesize $v_5$}
%\psfrag{v6}{\footnotesize $v_6$}
%\psfrag{v7}{\footnotesize $v_7$}
%\psfrag{v8}{\footnotesize $v_8$}
%\includegraphics*[scale=0.8]{bidirected_graph.eps}
%\end{center}
%\caption{Bidirected graph $\Sigma$} \label{fig:bigraph}
%\end{figure}
A \emph{walk} in a bidirected graph is a sequence  $(v_1, e_1, v_2, e_2,\ldots, e_{t-2}, v_{t-1}, e_{t-1}, 
v_t)$ where $v_i$ and $v_{i+1}$ are end-nodes of edge $e_i$ ($i=1,\ldots,t-1$), including the case where  
$v_i=v_{i+1}$ and $e_i$ is a half-edge. If $v_1=v_t$, then the walk is \emph{closed}. 
A walk which consists of only links and does not cross itself, that is $v_i\ne v_j$ for $i\ne j$, is a \emph{path}. A closed walk which does not cross itself (except at $v_1=v_t$) is called a \emph{cycle}. 
That is, a cycle can be a loop, a half-edge or a closed path. The \emph{sign of a cycle} is the product of 
the signs of its edges, so we have a \emph{positive cycle} if the number of negative edges in the cycle is 
even, otherwise the cycle is a \emph{negative cycle}. Obviously, a negative loop or a half-edge always 
makes a negative cycle. A bidirected graph is \emph{connected}, if there is a path between any two nodes. 
%A \emph{tree} is a connected graph which does not contain a cycle. 
%A connected graph containing exactly 
%one cycle is called \emph{1-tree}. This name is justified by the fact that a 1-tree consists of a tree and 
%one additional edge. If the unique cycle in a 1-tree is negative, then we will call it a 
%\emph{negative 1-tree}. 

The \emph{node-edge incidence matrix}  of a bidirected graph ${\Sigma}(V,E)$ 
is the $V\times{E}$ matrix $D_{\Sigma}$ with
\[
D_{\Sigma}(v,e)=
\left\{ 
\begin{array}{rl}
-1 & \textrm{if $v$ is an out-node of $e$}, \\
+1  & \textrm{if $v$ is an in-node of $e$}, \\
-2  & \textrm{if $e$ is a negative loop and $v$ is its out-node}, \\
+2  &  \textrm{if $e$ is a negative loop and $v$ is its in-node}, \\
0  & \textrm{otherwise,}
\end{array} \right.
\]
for any vertex $v\in V$ and any edge $e\in E$. 
% If $e$ is directed (i.e. positive) loop then we
% set  $A_{\Sigma}(v,e):=0$ for each vertex $v$.
%For example, the incidence matrix of the bidirected graph depicted in Figure~\ref{fig:bigraph} is: 
%\[
%A_{\Sigma}=\begin{tabular}{r|rrrrrrrrrrrrrr|} 
%\multicolumn{1}{c}{} & \multicolumn{1}{c}{$e_1$} &\multicolumn{1}{c}{$e_2$} & 
%\multicolumn{1}{c}{$e_3$} & \multicolumn{1}{c}{$e_4$} & \multicolumn{1}{c}{$e_5$}& 
%\multicolumn{1}{c}{$e_6$}& \multicolumn{1}{c}{$e_7$}& \multicolumn{1}{c}{$e_8$} & 
%\multicolumn{1}{c}{$e_9$} & \multicolumn{1}{c}{$e_{10}$} & \multicolumn{1}{c}{$e_{11}$} & 
%\multicolumn{1}{c}{$e_{12}$} & \multicolumn{1}{c}{$e_{13}$} & \multicolumn{1}{c}{$e_{14}$} \\ \cline{2-15} 
%$v_1$ &-1 & 0  &  1 & 0& 0 &  0 & 0  & 0 & 0  & 0 & 0 & 0 & 0   & 0 \\  
%$v_2$ & 1 & -1 &  1 & 0& 0 &  1 & 0  & 0 & 0  & 0 & 0 & 0 & 0   & 0 \\  
%$v_3$ & 0 & 1  &  0 & 1& -1 & 0 & 0  & 0 & 0  & 0 & 0 & 0 & 0   & 0 \\ 
%$v_4$ & 0 & 0  &  0 & 0& 1 &  0 & -1 & 1 & 0  & 0 & 0 & 0 & 0   & 0 \\ 
%$v_5$ & 0 & 0  &  0 & 0& 0 &  0 & 1  & 0 & 1  & 0 & 0 & 0 & 0   & 0 \\ 
%$v_6$ & 0 & 0  &  0 & 0& 0 &  0 & 0  & 1 & -1 & 1 & 0 & 0 & 0   & 0 \\ 
%$v_7$ & 0 & 0  &  0 & 0& 0 &  0 & 0  & 0 & 0  & 1 & 1 & 2 & 1   & 0 \\ 
%$v_8$ & 0 & 0  &  0 & 0& 0 &  0 & 0  & 0 & 0  & 0 & 1 & 0 & -1  & 1 \\ \cline{2-15} 
%\end{tabular}
%\]
The following operations are defined on bidirected graphs. \emph{Deletion of an edge} is simply the 
removal of the edge; \emph{deletion of a node} means that the node and all the edge-ends incident to it 
are removed. Thus incident half-edges or loops become loose edges, incident links become half-edges. 
Deletion of an edge or a node is equivalent to the deletion of the corresponding column or row from the 
node-edge incidence matrix. \emph{Switching} at a node is the operation when all the labels at the incident 
edge-ends are changed to the opposite. It corresponds to the multiplication by $-1$ of a row in the 
incidence matrix. Finally, \emph{contracting an edge} $e$ is the operation in which the end-nodes of $e$ 
are modified and $e$ is shrunk to zero length. For different types of edges 
contraction manifests itself differently. Specifically, if $e$ is a negative loop or a half-edge then the node 
incident to it is deleted together with all the edge-ends incident to it. If $e$ is a positive link, 
then its two end-nodes are identified and $e$ is deleted. If $e$ is a negative link, then first we switch 
at one of its end-nodes to make it a positive link, then contract it as defined for positive links. If $e$ is
a positive loop then $e$ is simply deleted. 
%
%There are more general graphs than bidirected graphs that can be discussed. One of these is the class of 
%signed graphs introduced by Harary~\cite{Harary:1953} which can be obtained from bidirected graphs by 
%ignoring the signs at the endpoints of each 
%edge and focusing on whether the edge is positive or negative. In other words, a bidirected graph can be 
%viewed as the directed version of a signed graph \cite{Zaslavsky:1991b}. More about signed graphs can be 
%found in~\cite{Zaslavsky:1998,Zaslavsky:1999}.
%
%%%%%%%%%%%%%%%%%%%%%%%%%%%%%%%%%%%%%%%%%%%%%%%%%%%%%%%%%%%%%%%%%%%%%%%%%%%%%%%%%%%%%%%%%%%%%%%%%%%%%%%%%%
%\subsection{Binet Matrices} \label{sec_binet}
Binet matrices are defined similarly as network matrices as follows:
\begin{definition}\label{def_binet}
Let $D_{\Sigma}$ be a full row rank node-edge incidence matrix of a bidirected graph $\Sigma$, $R$ be a 
basis of it and $D_{\Sigma}=[R|S]$. The matrix $B=R^{-1}S$ is called a \emph{binet matrix}.
\end{definition}
%The subgraph of $\Sigma$ induced by the edges corresponding to the columns of $R$ is called a 
%\emph{basis of the graph}, and will be denoted by $\Sigma(R)$. Its edges 
%are the \emph{basic edges} while the edges corresponding to columns of $S$ are the \emph{non-basic edges} 
%of $\Sigma$. The rows and columns of $B$ are identified with the basic and non-basic edges of $\Sigma$ 
%respectively. 
%$\Sigma(R)$ will always be a collection of unconnected negative 1-trees. The unconnected components are 
%called the \emph{basic components}, the unique cycles in the 1-trees are the \emph{basic cycles}. 
When in a bidirected graph $\Sigma$ the subgraph $\Sigma(R)$ is indicated for a basis $R$, we call it a 
\emph{binet representation} or a \emph{binet graph}.  
In Figure~\ref{fig:fundex} the binet graph for basic edges  $\{r_1, \ldots, r_8\}$ and non-basic edges  
$\{s_1,\ldots, s_6\}$ is shown, with the associated binet matrix. 
It is noted that in computing the entries of a binet matrix for a given basis, instead of using Definition~\ref{def_binet}
which involves the inverse of a matrix, there also exists a combinatorial algorithm described 
in~\cite{AppKot:2005,BolZasl:2005}. 

For a column $s$ of $S$, let $r_1, r_2, \ldots, r_t$ be the columns of $R$ for which the  corresponding 
component of vector $R^{-1}s$ is non-zero. The vectors $r_1, r_2, \ldots, r_t$ and $s$ form a minimal linearly
dependent set in $\mathbb{R}$. The subgraphs
of $\Sigma$ induced by sets of edges which correspond to minimally dependent sets of columns in $A_{\Sigma}$
have to be one of the following three types as shown in~\cite{Kotnyek:2002,Zaslavsky:1982}:
\begin{itemize}
\item[(i)] a positive cycle,
\item[(ii)] a graph consisting of two negative cycles which have exactly one common node,
\item[(iii)] a graph consisting of two node-disjoint negative cycles connected with a path which has no common node with 
the cycles except its end-nodes.
\end{itemize}
Graphs in categories (ii) and (iii) are called \emph{handcuffs} of type I and II respectively.
For example in the bidirected graph illustrated in Figure~\ref{fig:fundex} the subgraph induced by
the edges $\{r_1, r_3, r_4, s_1\}$ is a positive cycle, while the sets of edges $\{r_1, r_2, r_3, s_3\}$ and
$\{r_1,r_2, r_3, r_4, r_5, r_6, s_2\}$ induce handcuffs of type I and II respectively. 

Some results concerning binet matrices which will be useful  are the following. Proofs can be found in~\cite{AppKot:2005,Kotnyek:2002}.
\begin{theorem} Binet matrices are closed under the following operations: 
\smallskip

\noindent(a) Switching at a node of a binet graph. 

\noindent(b) Multiplying a row or column with $-1$. 

\noindent(c) Deleting a row or a column.

\noindent(d) Pivoting (in $\mathbb R$) on a nonzero element. 
\end{theorem}

Switching at a node does not change the matrix, the new binet graph represents the same matrix. Multiplying a row or 
column with $-1$ is equivalent to reversing the orientation of the corresponding basic or non-basic edge. Deleting a 
column is simply deleting the corresponding non-basic edge, while deleting a row amounts to contracting the 
corresponding basic edge. Pivoting on an element in row $r$ and column $s$ means that these edges are exchanged in the 
basis.

%%%%%%%%%%%%%%%%%%%%%%%%%%%%%%%%%%%%%%%%%%%%%%%%%%%%%%%%%%%%%%%%%%%%%%%%%%%%%%%%%%%%%%%%%%%%%%%%%%%%%%%%%%%%%%%%%%%%%%%%%%%%%%%%
\subsection{Decomposition of Totally Unimodular Matrices}   
A matrix $A$ is totally unimodular if each square submatrix of $A$ has determinant $0,+1,$ or $-1$. 
There are numerous other characterisations of the class of TU matrices 
(see \cite{NemWols:1988,Schrijver:98}). 
%The importance of TU matrices in combinatorial 
%optimisation stems mainly from 
%Hoffman and Kruskal's characterisation \cite{HofKru:56}, which states that an integral matrix $A$ is TU 
%if and only if for each integral vector $b$ the polyhedron $P=\{x|x\geq{0};Ax\leq{b}\}$ is integral.
%As a direct application of this, the integer programming problem $\max \{cx | x\in P, x \textrm { integer}\}$
%can be solved as a linear programming problem in polynomial time. 
%The most important subclass of totally unimodular matrices is that of network matrices 
%(Tutte \cite{Tutte:1965}), that is any network matrix is totally unimodular. 
The following decomposition theorem for TU  matrices proved by Seymour~\cite{Seymour:1980} plays a central
role in this work, and also yields a polynomial-time test for total unimodularity. 
\begin{theorem} \label{Seymour_matrix}
Any totally unimodular matrix is up to row and column permutations and scaling by $\pm{1}$ factors a network 
matrix, or the transpose of such a matrix, or the matrix $B_1$ or $B_2$ of (\ref{eq_B1}) and
(\ref{eq_B2}), or may be constructed 
recursively from these matrices using  matrix $1$-, $2$- and $3$-sums (see Definition~\ref{def_k-sums}). 
\end{theorem}  
Matrices $B_1$ and $B_2$ are binet matrices, as it is indicated by the corresponding binet graphs shown
in (\ref{eq_B1}) and (\ref{eq_B2}).
The above theorem is essentially a direct consequence of a decomposition theory for matroids associated
with TU  matrices, the so-called \emph{regular matroids}. Specifically Seymour characterised the class of
regular matroids by defining certain operations called $k$-sums,  
such that every regular matroid can be decomposed into a set of elementary building blocks via these
operations, if and only if these blocks satisfy certain properties. 
\begin{equation}\label{eq_B1}
%\hfill
B_1=\begin{tabular}{r|rrrrr|} 
\multicolumn{1}{c}{} &  \multicolumn{1}{c}{$s_1$} &\multicolumn{1}{c}{$s_2$} & 
\multicolumn{1}{c}{$s_3$} & \multicolumn{1}{c}{$s_4$} & \multicolumn{1}{c}{$s_5$}\\ \cline{2-6} 
$r_1$ & 1  & 0  & 0  & 1  & -1  \\ %\cline{2-6} 
$r_2$ & -1 & 1  & 0  & 0  & 1 \\ %\cline{2-6} 
$r_3$ & 1  & -1 & 1  & 0  & 0 \\ %\cline{2-6} 
$r_4$ & 0  & 1  & -1 & 1  & 0 \\ %\cline{2-6} 
$r_5$ & 0  & 0  & 1  & -1 & 1 \\ \cline{2-6} 
\end{tabular}\mspace{120mu}
\raisebox{-10ex}{\resizebox*{0.6\width}{!}{
\psfrag{r1}{\footnotesize $r_1$}
\psfrag{r2}{\footnotesize $r_2$}
\psfrag{r3}{\footnotesize $r_3$}
\psfrag{r4}{\footnotesize $r_4$}
\psfrag{r5}{\footnotesize $r_5$}
\psfrag{s1}{\footnotesize $s_1$}
\psfrag{s2}{\footnotesize $s_2$}
\psfrag{s3}{\footnotesize $s_3$}
\psfrag{s4}{\footnotesize $s_4$}
\psfrag{s5}{\footnotesize $s_5$}
\includegraphics*[scale=0.9]{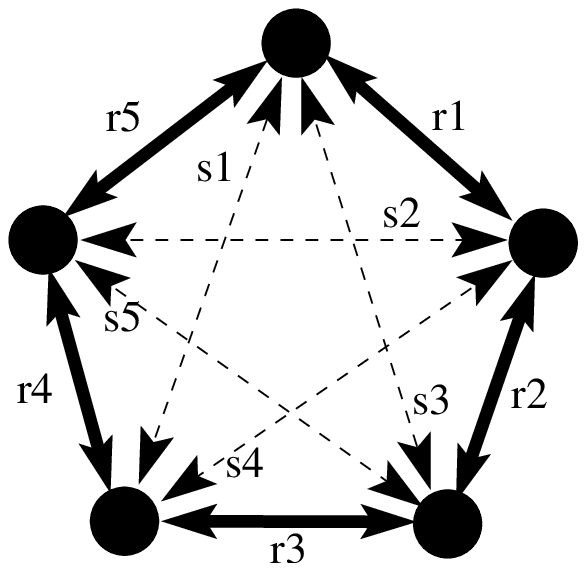}\hfill
}}
\hspace*{\fill}
\end{equation}
\begin{equation}\label{eq_B2}
%\hfill
B_2=\begin{tabular}{c|ccccc|} 
\multicolumn{1}{c}{} & \multicolumn{1}{c}{$s_1$} &\multicolumn{1}{c}{$s_2$} & 
\multicolumn{1}{c}{$s_3$} & \multicolumn{1}{c}{$s_4$} & \multicolumn{1}{c}{$s_5$}\\ \cline{2-6} 
$r_1$ & 1 & 1 & 1 & 1 & 1  \\ %\cline{2-6} 
$r_2$ & 1 & 1 & 1 & 0 & 0 \\ %\cline{2-6} 
$r_3$ & 1 & 0 & 1 & 1 & 0 \\ %\cline{2-6} 
$r_4$ & 1 & 0 & 0 & 1 & 1 \\ %\cline{2-6} 
$r_5$ & 1 & 1 & 0 & 0 & 1 \\ \cline{2-6} 
\end{tabular}\mspace{120mu}
\raisebox{-10ex}{\resizebox*{0.6\width}{!}{
\psfrag{r1}{\footnotesize $r_1$}
\psfrag{r2}{\footnotesize $r_2$}
\psfrag{r3}{\footnotesize $r_3$}
\psfrag{r4}{\footnotesize $r_4$}
\psfrag{r5}{\footnotesize $r_5$}
\psfrag{s1}{\footnotesize $s_1$}
\psfrag{s2}{\footnotesize $s_2$}
\psfrag{s3}{\footnotesize $s_3$}
\psfrag{s4}{\footnotesize $s_4$}
\psfrag{s5}{\footnotesize $s_5$}
\includegraphics*[scale=0.9]{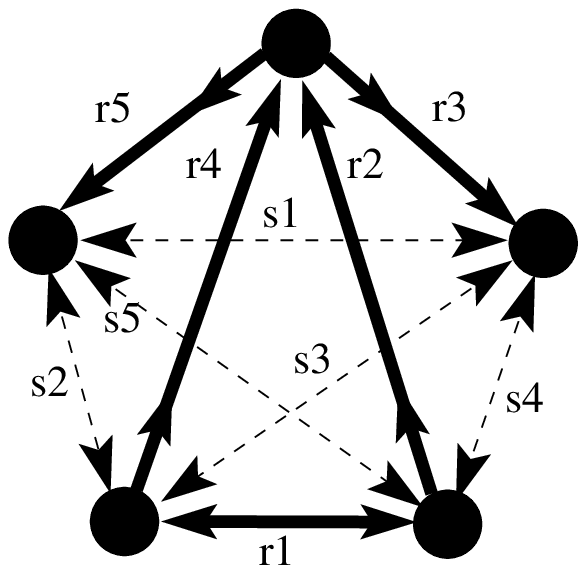}\hfill
}}
\hspace*{\fill}            
\end{equation}

In general, $k$-sum operations $(k=1,2,3)$ are defined in the more general theoretical framework of matroids, and 
here we basically treat the specialised version of this operation as applied to the compact representation 
matrices of regular matroids. Moreover, it can be shown that applying these operations on totally 
unimodular matrices preserves their total unimodularity.
\begin{definition}\label{def_k-sums}
If $A, B$ are matrices and $a,d$ and $b,c$ are column and row vectors of
appropriate size in $\mathbb R$ then 
\begin{description} \label{des_ksums}
\item[1-sum:] $A\oplus_1 B:=\begin{bmatrix}A & 0\\0& B\end{bmatrix}$
\item[2-sum:] $\begin{bmatrix}A & a\end{bmatrix}\oplus_2 \begin{bmatrix}b\\B\end{bmatrix}:=\begin{bmatrix}A & ab\\0& 
B\end{bmatrix}$
\item[3-sum:] $\begin{bmatrix}A & a & a\\c & 0 & 1\end{bmatrix}\oplus_3 \begin{bmatrix}1 & 0 & b\\d & d & B\end{bmatrix}:=
\begin{bmatrix}A & ab\\dc& B\end{bmatrix}$ or \\
\hspace*{.13in}$\begin{bmatrix}A & 0 \\b & 1 \\ c & 1 \end{bmatrix}\oplus^3 \begin{bmatrix}1 & 1 & 0\\
a& d & B\end{bmatrix}:=\begin{bmatrix}A & 0 \\ D & B\end{bmatrix}$ \\
where in the $\oplus^3$-sum row vectors $b$ and $c$ and column vectors $a$ and $d$ are submatrices of $D$ and the $2\times{2}$ matrix $\bar{D}$ is the intersection of 
rows $b$ and $c$ with columns $a$ and $d$. Further the rank of   $D=[a|d]\bar{D}^{-1}[\frac{b}{c}]$ is two. Note that there are two alternative definitions for $3$-sum, 
distinguished by $\oplus_3$ and $\oplus^3$. The indices of the isolated columns and rows in the 2-sum and 3-sum operations, will be called {\em connecting elements}.
\end{description}  
\end{definition}

The definition of the $k$-sums may seem complicated at first glance, but they essentially provide a 
way to decompose a TU matrix into smaller TU matrices given that the matrix admits such a decomposition. 
Specifically suppose that we have a TU matrix $N$ which under row and column permutations can take the form
\begin{equation}\label{eq_example_matrix}
N = \begin{bmatrix}A & D_1 \\ D_2 & B\end{bmatrix}
\end{equation}
and the following two conditions are satisfied:
\begin{itemize}
\item[(i)] number of rows and columns of both $A$ and $B$ $> k$,
\item[(ii)] $\rank(D_1) + \rank(D_2) = k-1$ where $D_1,D_2$ are viewed over $GF(2)$. 
\end{itemize}
Then the matrix $N$ of (\ref{eq_example_matrix}) can be decomposed under a $k$-sum operation into two 
matrices of smaller size which are submatrices of $N$, 
preserving total unimodularity. In the case of $3$-sum we note from the definition that there
are two alternative operations,  reflecting the fact that condition (ii) above can be satisfied
in two different ways (i.e. $\rank(D_1)=\rank(D_2) = 1$ or $\rank(D_1)=0, \rank(D_2) = 2$). However
it can be shown that when the matrices are TU both definitions of $3$-sum are equivalent under 
pivoting in either $GF(2)$ or $\mathbb{R}$. 
(The regular matroid decomposition theorem of Seymour, $k$-sums of matrices and their corresponding matroids, 
and decomposition theory for matroids in general is treated extensively in ~\cite{Oxley:92,Truemper:98}.)

%%%%%%%%%%%%%%%%%%%%%%%%%%%%%%%%%%%%%%%%%%%%%%%%%%%%%%%%%%%%%%%%%%%%%%%%%%%%%%%%%%%%%%%%%%
\section{$k$-sum of Network and Binet Matrices}\label{sec_ksums}
In this section we will examine the operation of $k$-sums of matrices, both network and
binet. We will show whether the resulting matrix
is a network or binet matrix, or does not belong to either class. 
Algebraic proofs as well as graphical representations of the associated 
operations on these matrices are presented.

%%%%%%%%%%%%%%%%%%%%%%%%%%%%%%%%%%%%%%%%%%%%%%%%%%%%%%%%%%%%%%%%%%%%%%%%%%%%%%%%%%%%%%%%%%%%%%%%%%%%%%%%%%%%%%%%%%%%%%%%%%%%%%%%
\subsection{$k$-sums of Network Matrices}
Here it is proved that network matrices are closed under the $k$-sum operations. Since network 
matrices are the compact representation matrices of graphic matroids, a direct consequence of these 
results is the well known fact (see \cite{Oxley:92}) that graphic matroids are closed under $k$-sums. 
However the analytical methodology in the  proofs that will be given here, will be used in the sections that
will follow where the binet, and the more general tour matrix case is treated. Moreover since the proof
is constructive, it is used in the algorithm for composing the bidirected graph of a TU matrix which
will be presented in section~\ref{subsubsec_algorithm}.

%%%%%%%%%%%%%%%%%%%%%%%%%%%%%%%%%%%%%%%%%%%%%%%%%%%%%%%%%%%%%%%%%%%%%%%%%%%%%%%%%%%%%%%%%%%%%%%%%%%%%%%%%%%%%%%%%%%%%%%%%%%%%%%%
\subsubsection{Network $\oplus_3$ Network} \label{sec_net3net}
The most general case of 3-sum will be examined since the other sum operations follow. 
\begin{lemma}\label{lem_net-3sum-net}
If $N_1$, $N_2$ are network matrices such that 
\[
N_1=
\raisebox{5pt}{$
\begin{array}{r}
\vspace{3.2mm}    \\
e_3
\end{array}
\hspace{-2.5mm}
\begin{array}{c}
\begin{array}{crr}
\hspace{1mm} & e_1 & \hspace{-1.7mm} e_2 
\end{array} \\
\left[ \begin{array}{ccc}
A & a & a\\
c & 0 & 1
\end{array}
   \right],
\end{array} 
$}
N_2 = 
\raisebox{5pt}{$
\begin{array}{r}
\vspace{-5mm}
\hspace{-3.2mm}    \\
f_3
\end{array}
\hspace{-2.5mm}
\begin{array}{c}
\begin{array}{crr}
\hspace{-3.2mm} f_1 & \hspace{-1.3mm} f_2 &  
\end{array} \\
\left[ \begin{array}{ccc}
1 & 0 & b\\
d & d & B
\end{array}
   \right],
\end{array} 
$}\]
then $N=N_1 \oplus_3 N_{2}$ is a network matrix. 
\end{lemma}
\begin{proof}
Because of the definition of the $3$-sum operation we have that in a possible graphical representation of $N_1$ the 
fundamental cycle of $e_1$ consists of the edges that correspond to non-zero elements in $a$. The fundamental cycle of $e_2$ has all 
these edges and $e_3$. This means that $e_1$, $e_2$ and $e_3$ should form a triangle. Similarly, $f_1$, $f_2$ and 
$f_3$ form a triangle in any network representation of $N_2$. Let now $[R_1|S_1]$ and $[R_2|S_2]$ be the incidence matrices 
associated with $N_1$ and $N_2$, respectively, where 
after permutations and/or multiplications of rows with $\pm{1}$ we can write:
\begin{equation} \label{eq_fourtythree}
[R_1|S_1]=
\begin{array}{c}
\begin{array} {crrrc}
& \hspace{10mm} e_3 & \hspace{0.37in} e_1 & \hspace{1mm} e_2 &
\end{array}\\
\left[
\begin{array}{cr|crr}
r_1        & -1 & s_1   & 0   & -1\\
r_1'       & 1  & s_1'  & -1  & 0\\
r_1''      & 0  & s_1'' & 1   & 1\\
{R_1}'  & \mathbf{0} &{S_1}'   & \mathbf{0} & \mathbf{0} 
\end{array} \right]
\end{array}, \quad
[R_2|S_2]=
\begin{array}{c}
\begin{array} {ccccc}
\hspace{-2mm} f_3  &  & \hspace{5.5mm}  f_1 & \hspace{1mm} f_2 &
\end{array}\\
\left[
\begin{array}{rc|rrc}
0      & r_2   & -1  & -1 & s_2\\
-1     & r_2'  &  0  & 1  & s_2'\\
1      & r_2'' &  1  & 0  & s_2''\\
\mathbf{0}  & {R_2}'   & \mathbf{0} & \mathbf{0} & {S_2}'
\end{array} \right]
\end{array}
\end{equation}
where $\mathbf{0}$ is a vector or matrix of zeros of appropriate size, $r_{i},r_{i}',r_{i}'',s_{i},s_{i}'$ and $s_{i}''$ are row vectors and ${R_{i}}', {S_{i}}'$ are  matrices 
of appropriate size $(i=1,2)$. By the definition of network matrices the following two equations hold:
\begin{equation} \label{eq_fourtyfive}
R_1N_1=S_1, \qquad
R_2N_2=S_2
\end{equation}
For $N_1$ using (\ref{eq_fourtythree}) and (\ref{eq_fourtyfive}) we have that:
\begin{equation*}
\left[
\begin{array}{c|r}
r_1  &  -1\\ 
r_1' &   1\\
r_1''&   0\\ \hline
{R_1}' & \mathbf{0}
\end{array} \right]
\left[
\begin{array}{c|c|c}
A  &  a  & a \\ \hline
c  &  0  & 1
\end{array}\right]=
\left[ \begin{array}{c|r|r}
s_1       &   0  & -1\\
s_1'      &  -1  &  0\\
s_1''     &   1  &  1\\ \hline
{S_1}' & \mathbf{0} & \mathbf{0}
\end{array} \right]
\end{equation*}
where upon decomposing the block matrix multiplications we derive the following equations.
\begin{align} \label{eq_fourtyeight}
\left[ \begin{array}{c}
r_1\\
r_1'\\
r_1''
\end{array} \right]A+
\left[ \begin{array}{r}
-1\\
1\\
0
\end{array} \right]c=
\left[ \begin{array}{c}
s_1\\
s_1'\\
s_1''
\end{array} \right],\quad
\left[ \begin{array}{c}
r_1\\
r_1'\\
r_1''
\end{array} \right]a =
\left[ \begin{array}{r}
0\\
-1\\
1
\end{array} \right]\\ \nonumber  
\left[ \begin{array}{c}
r_1\\
r_1'\\
r_1''\\
\end{array} \right]a+
\left[\begin{array}{r}
-1\\
1\\
0
\end{array} \right]=
\left[ \begin{array}{r}
-1\\
0\\
1
\end{array} \right], \quad
{R_1}'A={S_1}', \quad
{R_1}'a=\mathbf{0}
\end{align}
Similarly, for $N_2$ using (\ref{eq_fourtythree}) and (\ref{eq_fourtyfive}) we have
\begin{equation*}
\left[ \begin{array}{r|c}
0    &   r_2\\
-1   &   r_2'\\
1    &   r_2''\\ \hline
\mathbf{0} & {R_2}'
\end{array} \right]
\left[ \begin{array}{c|c|c}
1    &  0   &  b\\ \hline
d    &  d   &  B
\end{array} \right]=
\left[ \begin{array}{r|r|c}
-1   &  -1  & s_2\\
0    &   1  & s_2'\\
1    &   0  & s_2''\\ \hline
\mathbf{0} & \mathbf{0} & {S_2}'
\end{array} \right]
\end{equation*}
so that
\begin{align} \label{eq_fiftyfour}
\left[ \begin{array}{r}
0\\
-1\\
1
\end{array} \right]+
\left[ \begin{array}{c}
r_2\\
r_2'\\
r_2''
\end{array} \right]d=
\left[ \begin{array}{r}
-1\\
0\\
1
\end{array} \right], \quad
\left[ \begin{array}{c}
r_2\\
r_2'\\
r_2''
\end{array} \right]d=
\left[ \begin{array}{r}
-1\\
1\\
0
\end{array} \right] \quad \\ \nonumber
\left[ \begin{array}{r}
0\\
-1\\
1
\end{array} \right]b+
\left[ \begin{array}{c}
r_2\\
r_2'\\
r_2''
\end{array} \right]B=
\left[ \begin{array}{c}
s_2\\
s_2'\\
s_2''
\end{array} \right],\quad
{R_2}'d=\mathbf{0}, \quad
{R_2}'B={S_2}'
\end{align}
Using block matrix multiplication and equations in (\ref{eq_fourtyeight}) and (\ref{eq_fiftyfour}), it is easy to show that the following equality holds:
\begin{equation} \label{eq_fiftynine}
\underbrace{
\left[ \begin{array}{c|c}
r_1  & r_2\\
r_1' & r_2'\\
r_1''& r_2''\\ \hline
{R_1}' & \mathbf{0}\\ \hline
\mathbf{0} & {R_2}'
\end{array} \right]
}_{R'}\underbrace{
\left[ \begin{array}{c|c}
A  & ab\\ \hline
dc & B
\end{array} \right]
}_{N}=
\underbrace{
\left[ \begin{array}{c|c}
s_1  & s_2\\
s_1' & s_2'\\
s_2''& s_2''\\ \hline
{S_1}' & \mathbf{0}\\ \hline
\mathbf{0} & {S_2}'
\end{array} \right]
}_{S'}\end{equation}
The matrix $[R'|S']$ is the incidence matrix of a directed graph since each column contains a $+1$ and a $-1$. It remains to be shown that the matrix $\hat{R}$ obtained by deleting one 
row of $R'$ is non-singular. If we delete the first row of $R'$ we have that:
\begin{equation*}
\hat{R}=
\left[ \begin{array}{c|c}
r_1'  & r_2'\\
r_1'' & r_2''\\
{R_1}' & \mathbf{0}\\ \hline
\mathbf{0} & {R_2}'
\end{array} \right]
\end{equation*}
If we delete the first row from $R_1$ then we obtain the matrix
$
\left[ \begin{array}{cc}
r_1'   & 1\\
r_1''   & 0\\
{R_1}' & \mathbf{0}
\end{array} \right]
$
which is a non-singular one.
Expanding now the determinant of that matrix along the last column we can see that the matrix 
$\left[ \begin{array}{c}
r_1''\\
{R_1}'
\end{array} \right]
$ 
is also non-singular.
Therefore, within the
submatrix $\left[ \begin{array}{c}
r_1'\\
r_1''\\
{R_1}'
\end{array} \right]$ of  $\hat{R}$, $r_1'$ can be written as a linear combination of the other rows:
\begin{equation} \label{eq_ninetyfive}
r_1'+u~r_1''+q {R_1}'=0
\end{equation}
where $u$ is a scalar, and $q$ is a column vector of appropriate size with elements in $\mathbb{R}$. Also, we have that $u\neq{0}$ since if we delete $e_3$ in 
$R_1$ then the matrix obtained corresponds to a forest in which the nodes which correspond to rows $r_1'$ and $r_1''$ belong to the same tree of that forest. We 
denote the determinant of $\hat{R}$ by $det[\hat{R}]$. Using (\ref{eq_ninetyfive}) we get:
\begin{equation} \label{eq_ninetysix}
det[\hat{R}]=
det\left[\begin{array}{c|c}
r_1'  & r_2'\\
r_1'' & r_2''\\
{R_1}' & \mathbf{0}\\ \hline
\mathbf{0} & {R_2}'
\end{array} \right]=
det\left[\begin{array}{c|c}
0     & r_2'+u~r_2''\\
r_1'' & r_2''\\
{R_1}' & \mathbf{0}\\ \hline
\mathbf{0} & {R_2}'
\end{array} \right]  = 
det\left[\begin{array}{c|c}
r_1'' & r_2''\\
{R_1}' & \mathbf{0}\\ \hline
0     & r_2'+u~r_2''\\
\mathbf{0} & {R_2}'
\end{array} \right]  
\end{equation}
So, matrix $\hat{R}$ is block diagonal and its blocks are square. Thus:
\begin{equation} \label{eq_ninetyseven}
det\left[ \hat{R}\right]=
det\left[ \begin{array}{c}
r_1''\\
{R_1}'
\end{array} \right]
det\left[ \begin{array}{c}
r_2'+u~r_2''\\
{R_2}'
\end{array} \right]=
det\left[ \begin{array}{c}
r_1''\\
{R_1}'
\end{array} \right]
\left(det\left[ \begin{array}{c}
r_2'\\
{R_2}'
\end{array} \right]+
u\,det\left[ \begin{array}{c}
r_2''\\
{R_2}'
\end{array}\right]
\right)
\end{equation}
If we delete from $R_2$ its first row then the matrix so obtained is non-singular and, since it is a submatrix 
of a TU matrix, it has to be TU as well, i.e. its determinant should be equal to $\pm{1}$. Expanding the determinant of that 
matrix along its first column we take:
\begin{equation} \label{eq_ninetyeight}
det \left[ \begin{array}{c}
r_2'\\
{R_2}'
\end{array} \right]+
det\left[ \begin{array}{c}
r_2''\\
{R_2}'
\end{array} \right]=\pm{1}
\end{equation}
Furthermore 
$det\left[ \begin{array}{c}
r_2'\\
{R_2}'
\end{array} \right],
det\left[ \begin{array}{c}
r_2''\\
{R_2}'
\end{array} \right]\in{\{0,\pm{1}\}}$
since the corresponding matrices are TU. 
From (\ref{eq_ninetyeight}) we see that exactly one of these matrices has a nonzero determinant. Combining this with 
(\ref{eq_ninetyseven}) and the fact
that $u\neq{0}$ we have that $\hat{R}$ is nonsingular.\\
Finally, it is obvious that the matrix $[R'|S']$ contains a $-1$ and a $+1$ in each column since its columns are 
columns of $[R_1|S_1]$ and $[R_2|S_2]$. We can conclude that the 3-sum of two network matrices is a 
network matrix with incidence matrix $[R'|S']$.
\end{proof}
%
%\noindent
%\textbf{Graphical Representation of Network $\oplus_{3}$ Network:}\\
%\noindent
%Using $[R'|S']$ from (\ref{eq_fiftynine}), we can draw a network representation of $N$ using the network representations 
%associated with $N_1$ and $N_2$. Gluing together the triangles $(e_1,e_2,e_3)$ and $(f_1,f_2,f_3)$ such that $e_3$ 
%meets $f_2$, $e_1$ meets $f_3$ and $e_2$ meets $f_1$ is the procedure that gives rise to $[R'|S']$ which is described 
%in (\ref{eq_fiftynine}). An illustrative example is given in Figure \ref{fig:net3net}, where the solid edges correspond to basic edges. 
%\begin{figure}
%\begin{center}
%\centering
%\psfrag{sum}{\footnotesize $\oplus_{3}$}
%\psfrag{=}{\footnotesize $=$}
%\psfrag{e1}{\footnotesize $e_1$}
%\psfrag{e2}{\footnotesize $e_2$}
%\psfrag{e3}{\footnotesize $e_3$}
%\psfrag{f1}{\footnotesize $f_1$}
%\psfrag{f2}{\footnotesize $f_2$}
%\psfrag{f3}{\footnotesize $f_3$}
%\psfrag{r}{\footnotesize $r$}
%\psfrag{s}{\footnotesize $s$}
%\psfrag{a}{\footnotesize $a$}
%\psfrag{d}{\footnotesize $d$}
%\includegraphics*[scale=0.20]{net_3sum_net.eps}
%\end{center}
%\caption{The network representation of the 3-sum of two network matrices.\label{fig:net3net}}
%\end{figure}

\begin{theorem} \label{th_net_cl}
Network matrices are closed under $k$-sums $(k=1,2,3)$.
\end{theorem}
\begin{proof}
For $k=1$ it is straightforward. For $k=2$ it is enough to observe that if 
$N_1=\begin{bmatrix}A & a\end{bmatrix}, N_2=\begin{bmatrix}b\\B\end{bmatrix}$ are network matrices, then the 
matrices $\bar{N}_1=\begin{bmatrix}A & a & a\\0 & 0 & 1\end{bmatrix}$ and $\bar{N}_2=\begin{bmatrix}1 & 0 & b\\0 & 0 & B\end{bmatrix}$ are network matrices too, 
since we have only duplicated columns and added unitary rows and columns. But then $N_1\oplus_2{N_2}=\bar{N}_1\oplus_3{\bar{N}_2}$ which we know from 
Lemma~\ref{lem_net-3sum-net} to be network. For the alternative $3$-sum operation, since network matrices are closed under pivoting the result follows.
\end{proof}

%%%%%%%%%%%%%%%%%%%%%%%%%%%%%%%%%%%%%%%%%%%%%%%%%%%%%%%%%%%%%%%%%%%%%%%%%%%%%%%%%%%%%%%%%%%%%%%%%%%%%%%%%%%%%%%%%%%%%%%%%%%%%%%%
\subsection{$k$-sums of Network and Binet Matrices}
In this section we examine the $k$-sums between network and binet matrices. We prove that the result is always a 
binet matrix and we provide the associated bidirected graph representations. 

%%%%%%%%%%%%%%%%%%%%%%%%%%%%%%%%%%%%%%%%%%%%%%%%%%%%%%%%%%%%%%%%%%%%%%%%%%%%%%%%%%%%%%%%%%%%%%%%%%%%%%%%%
\subsubsection{Network $\oplus_3$ Binet}
Let's assume that $N_2$ of Lemma \ref{lem_net-3sum-net} is a binet matrix instead of a network matrix; then in a possible 
representation of it, its edges could be not 
only links but also loops and half edges. Most importantly, because of the structure of matrix $N_2$ we have that the edges 
$f_1$, $f_2$ and $f_3$ should be of a specific
type (loop, link, or half-edge) in order to form a binet representation of $N_2$. We examine below all the possible cases.

If $f_3$ is a link in the cycle (and then we can assume that it is a positive link), then $f_1$ and $f_2$ cannot be 
half-edges, because the fundamental circuit of a half-edge uses all the cycle edges, and the values on the cycle edges 
determined by the fundamental circuit are halves, so there can be neither 0 nor 1 in the row $f_3$ and columns $f_1$ and 
$f_2$ of $N_2$. Furthermore, $f_2$ cannot be a loop, because the fundamental circuit of any loop uses all cycle edges, 
despite the $0$ in the corresponding position of the matrix. So either both $f_1$ and $f_2$ are links, or $f_1$ is a loop 
and $f_2$ is a link. If they are both links, then they are both positive or both negative. Otherwise the fundamental 
circuit of one of them would use the negative edge in the cycle, the other would not, and they use the same edges except 
for the positive $f_3$. Moreover, $f_1,f_2$ and $f_3$ must form a triangle, so by a switching at a node we can make 
both $f_1$ and $f_2$ positive.

If $f_3$ is a loop, then $f_1$ cannot be a half-edge, because then the entry in row $f_3$ and column $f_1$ of $N_2$ 
would be a half. If $f_1$ is a loop, then vector $d$ of $N_2$ contains $\pm 2$ entries, but this is impossible because then 
$f_2$ would be an edge whose fundamental
circuit uses non-cycle edges twice but does not use the basic cycle (which is $f_3$). 
So $f_1$ must be a link, which implies that $f_2$ is also a link, and  $f_1$ is negative and $f_2$ is positive, because 
the fundamental circuit of $f_1$ uses the basic cycle, that of $f_2$ does not.

If $f_3$ is a half-edge, then $f_2$ must be a positive link, as its fundamental circuit does not use the basic cycle formed by
$f_3$. This also implies that $f_1$ is a half-edge.

If $f_3$ is a non-basic link, then $f_1$ cannot be a loop, as then it would have $\pm{2}$ on $f_3$ in the fundamental 
circuit. So either $f_1$ is a link and then $f_2$ is a link or a loop; or $f_1$ is a half-edge in which case $f_2$ is 
also a half-edge.

\begin{figure}[h]
\begin{center}
\centering
\psfrag{sum}{\footnotesize $\oplus_{3}$}
\psfrag{=}{\footnotesize $=$}
\psfrag{e1}{\footnotesize $e_1$}
\psfrag{e2}{\footnotesize $e_2$}
\psfrag{e3}{\footnotesize $e_3$}
\psfrag{f1}{\footnotesize $f_1$}
\psfrag{f2}{\footnotesize $f_2$}
\psfrag{f3}{\footnotesize $f_3$}
\psfrag{r}{\footnotesize $r$}
\psfrag{s}{\footnotesize $s$}
\psfrag{a}{\footnotesize $a$}
\psfrag{d}{\footnotesize $d$}
\includegraphics*[scale=0.20]{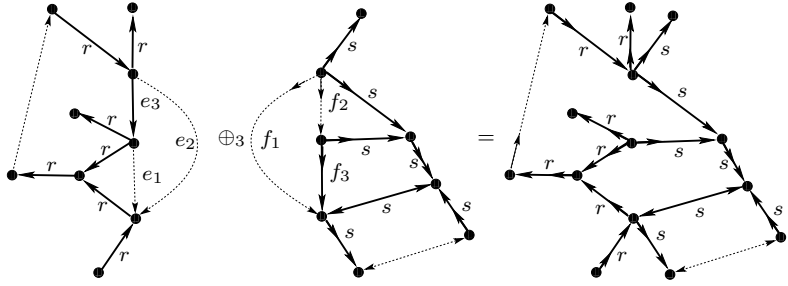}
\end{center}
\caption{The binet representation of the 3-sum of a network and a binet matrix. The case when $f_1,f_2, f_3$ are 
links.\label{fig:net3bin1}}
\end{figure}

\begin{figure}[h]
\begin{center}
\centering
\psfrag{sum}{\footnotesize $\oplus_{3}$}
\psfrag{=}{\footnotesize $=$}
\psfrag{e1}{\footnotesize $e_1$}
\psfrag{e2}{\footnotesize $e_2$}
\psfrag{e3}{\footnotesize $e_3$}
\psfrag{f1}{\footnotesize $f_1$}
\psfrag{f2}{\footnotesize $f_2$}
\psfrag{f3}{\footnotesize $f_3$}
\psfrag{r}{\footnotesize $r$}
\psfrag{s}{\footnotesize $s$}
\psfrag{a}{\footnotesize $a$}
\psfrag{d}{\footnotesize $d$}
\includegraphics*[scale=0.20]{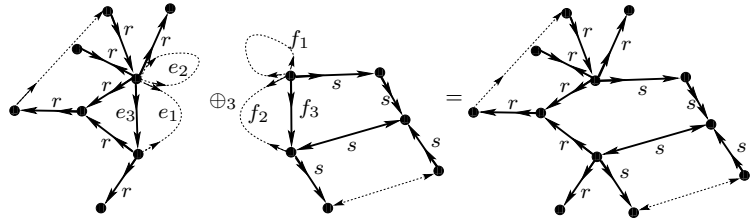}
\end{center}
\caption{The binet representation of the 3-sum of a network and a binet matrix.The case when $f_1$ is a loop, 
$f_2$ is a negative link,  $f_3$ is a positive link.\label{fig:net3bin2}}
\end{figure}

Therefore the cases that may appear are the following six: 
\begin{itemize}
\item[(a)] $f_3$ is a positive link in the cycle and $f_1$, $f_2$ are positive links; 
\item[(b)] $f_3$ is a positive link in the cycle, $f_1$ is a negative loop and $f_2$ is a negative link; 
\item[(c)] $f_3$ is a negative loop, $f_1$ is a negative link and $f_2$ is a positive link; 
\item[(d)] $f_1$, $f_3$ are half-edges and $f_2$ is a positive link; 
\item[(e)] $f_3$ is a non-cycle link, $f_1$ is a link and $f_2$ is a link or a negative loop; and 
\item[(f)] $f_3$ is a non-cycle link and $f_1$, $f_2$ are half-edges.
\end{itemize}          
\begin{lemma}\label{lem_net-3sum-binet}
If $N_1$ is a network matrix and $N_2$ is a binet matrix such that 
\[
N_1=
\raisebox{5pt}{$
\begin{array}{r}
\vspace{3.2mm}    \\
e_3
\end{array}
\hspace{-2.5mm}
\begin{array}{c}
\begin{array}{crr}
\hspace{1mm} & e_1 & \hspace{-1.7mm} e_2 
\end{array} \\
\left[ \begin{array}{ccc}
A & a & a\\
c & 0 & 1
\end{array}
   \right],
\end{array} 
$}
N_2 = 
\raisebox{5pt}{$
\begin{array}{r}
\vspace{-5mm}
\hspace{-3.2mm}    \\
f_3
\end{array}
\hspace{-2.5mm}
\begin{array}{c}
\begin{array}{crr}
\hspace{-3.2mm} f_1 & \hspace{-1.3mm} f_2 &  
\end{array} \\
\left[ \begin{array}{ccc}
1 & 0 & b\\
d & d & B
\end{array}
   \right],
\end{array} 
$}\]
then $N= N_1 \oplus_3 N_{2}$ is a binet matrix. 
\end{lemma}
\begin{proof}
Since $N_1$ is a network matrix we have that $e_1$, $e_2$ and $e_3$ should form a triangle. Therefore, w.l.o.g. we can assume for all the cases that the incidence matrix 
associated with the network matrix $N_1$ is the following one:
\begin{equation} \label{eq_st1} 
[R_1|S_1]=
\begin{array}{c}
\begin{array} {crrrc}
& \hspace{10mm} e_3 & \hspace{0.37in} e_1 & \hspace{1mm} e_2 &
\end{array}\\
\left[
\begin{array}{cr|crr}
r_1        & -1 & s_1   & 0   & -1\\
r_1'       & 1  & s_1'  & -1  & 0\\
r_1''      & 0  & s_1'' & 1   & 1\\
{R_{1}}'  & \mathbf{0} &{S_{1}}'   & \mathbf{0} & \mathbf{0} 
\end{array} \right],
\end{array}
\end{equation} where $\mathbf{0}$ is a zero matrix, $r_{i},r_{i}',r_{i}'',s_{i},s_{i}'$ and $s_{i}''$ are vectors and $R_{i}'$ 
and $S_{i}'$ are matrices of appropriate size $(i=1,2)$.

\noindent
{\bf Case (a):} For case (a) we have that the incidence matrix associated with the binet matrix $N_2$ can have the following form: 
\begin{equation*} 
[R_2|S_2]=
\begin{array}{c}
\begin{array} {ccccc}
\hspace{-2mm} f_3  &  & \hspace{5.5mm}  f_1 & \hspace{1mm} f_2 &
\end{array}\\
\left[
\begin{array}{rc|rrc}
0      & r_2   & -1  & -1 & s_2\\
-1     & r_2'  &  0  & 1  & s_2'\\
1      & r_2'' &  1  & 0  & s_2''\\
\mathbf{0}  & {R_{2}}'   & \mathbf{0} & \mathbf{0} & {S_{2}}'
\end{array} \right]
\end{array}
\end{equation*}
The proof for this case is very similar to the one regarding the $3$-sum of two network matrices in Lemma~\ref{lem_net-3sum-net}. 
Because of the structure of
matrix $N_2$, we have that $f_1$, $f_2$, and $f_3$ should form a triangle in any binet representation of $N_2$. Although we omit the 
full proof for this case because 
of its similarity to the one of Lemma \ref{lem_net-3sum-net}, we provide the incidence matrix matrix $[R'|S']$ of the binet graph 
associated with the binet matrix $N$ produced by the $3$-sum:
\begin{equation} \label{eq_sixty0}
[R'|S']=
\left[ \begin{array}{cc|cc}
r_1  & r_2  & s_1  & s_2\\
r_1' & r_2' & s_1' & s_2'\\
r_1''& r_2''& s_1''& s_2''\\ 
{R_{1}}' & \mathbf{0} & {S_{1}}'  & \mathbf{0}\\
\mathbf{0} & {R_{2}}' & \mathbf{0} & {S_{2}}'
\end{array} \right]
\end{equation}

\noindent
{\bf Case (b):} For this case we have that the incidence matrix associated with the binet matrix $N_2$ can have the following form:
\begin{equation} \label{eq_sixtyone}
[R_2|S_2]=
\begin{array}{c}
\begin{array} {ccccc}
\hspace{-2mm} f_3  &  & \hspace{5.5mm}  f_1 & \hspace{1mm} f_2 &
\end{array}\\
\left[ \begin{array}{rc|rrc}
-1 & r_2 & -2 & -1 & s_2\\
1  & r_2'& 0  & -1 & s_2'\\
\mathbf{0} & {R_{2}}' & \mathbf{0} & \mathbf{0} & {S_{2}}'
\end{array} \right]
\end{array}
\end{equation}
Initially, we convert the network representation $[{R_1}|{S_1}]$ of $N_1$ to a binet representation in 
which $e_2$ is a loop. This can be done so by introducing an artificial link 
parallel to $e_2$ and then contracting it. Thus, $e_1$ becomes a negative link, as contraction involves switching at the node to 
which $e_1$ and $e_2$ are incident. Graphically this case is illustrated in Figure~\ref{fig:net3bin2} which shows such an 
alternative binet representation of the matrix represented by the directed graph in Figure \ref{fig:net3bin1}. Therefore, the 
incidence matrix $[R_1|S_1]$ of the binet graph associated with $N_1$ can have the following form:   
\begin{equation} \label{eq_sixty}
[R_1|S_1]=
\begin{array}{c}
\begin{array} {crrrc}
& \hspace{10mm} e_3 & \hspace{0.37in} e_1 & \hspace{1mm} e_2 &
\end{array}\\
\left[ \begin{array}{cr|crr}
r_1  & -1  & s_1  & -1  & -2\\
r_1' &  1  & s_1' & -1  &  0\\
{R_{1}}' & \mathbf{0} & {S_{1}}' & \mathbf{0} & \mathbf{0}
\end{array} \right]
\end{array}
\end{equation}
We have that the following equations hold:
\begin{equation} \label{eq_sixtytwo}
R_1N_1=S_1, \qquad
R_2N_2=S_2
\end{equation}
From (\ref{eq_sixty}) and (\ref{eq_sixtytwo}) we have that:
\begin{equation*}
\left[
\begin{array}{c|r}
r_1  &  -1\\ 
r_1' &   1\\ \hline
{R_{1}}' & \mathbf{0}
\end{array} \right]
\left[
\begin{array}{c|c|c}
A  &  a  & a \\ \hline
c  &  0  & 1
\end{array}\right]=
\left[ \begin{array}{c|r|r}
s_1       &   -1  & -2\\
s_1'      &   -1  &  0\\\hline
{S_{1}}' & \mathbf{0} & \mathbf{0}
\end{array} \right],
\end{equation*}
where upon decomposing the block matrix multiplications we derive the following equations.    
\begin{gather} \label{eq_sixtyfive}
\left[ \begin{array}{c}
r_1\\
r_1'
\end{array} \right]A+
\left[ \begin{array}{r}
-1\\
1
\end{array} \right]c=
\left[ \begin{array}{c}
s_1\\
s_1'
\end{array} \right], \quad
\left[ \begin{array}{c}
r_1\\
r_1'
\end{array} \right]a =
\left[ \begin{array}{r}
-1\\
-1
\end{array} \right], \quad \nonumber \\
\left[ \begin{array}{c}
r_1\\
r_1'
\end{array} \right]a+
\left[\begin{array}{r}
-1\\
1
\end{array} \right]=
\left[ \begin{array}{r}
-2\\
0
\end{array} \right], \quad
{R_{1}}'A={S_{1}}', \quad
{R_{1}}'a=\mathbf{0}
\end{gather}
From (\ref{eq_sixtyone}) and (\ref{eq_sixtytwo}) we have that:
\begin{equation*} 
\left[ \begin{array}{r|c}
-2    &   r_2\\
1     &   r_2'\\ \hline
\mathbf{0} & {R_{2}}'
\end{array} \right]
\left[ \begin{array}{c|c|c}
1    &  0   &  b\\ \hline
d    &  d   &  B
\end{array} \right]=
\left[ \begin{array}{r|r|c}
-2   &  -1  & s_2\\
0    &   -1 & s_2'\\ \hline
\mathbf{0} & \mathbf{0} & {S_{2}}'
\end{array} \right]
\end{equation*}
and
\begin{gather} \label{eq_seventyone}
\left[ \begin{array}{r}
-1\\
1
\end{array} \right]+
\left[ \begin{array}{c}
r_2\\
r_2'
\end{array} \right]d=
\left[ \begin{array}{r}
-2\\
0
\end{array} \right], \quad
\left[ \begin{array}{c}
r_2\\
r_2'
\end{array} \right]d=
\left[ \begin{array}{r}
-1\\
-1
\end{array} \right], \quad \nonumber \\
\left[ \begin{array}{r}
-1\\
-1
\end{array} \right]b+
\left[ \begin{array}{c}
r_2\\
r_2'
\end{array} \right]B=
\left[ \begin{array}{c}
s_2\\
s_2'
\end{array} \right], \quad
{R_{2}}'d=\mathbf{0}, \quad
{R_{2}}'B={S_{2}}'
\end{gather}
Using block matrix multiplication and the equations in (18) and (19), the following equality holds:
\begin{equation*} 
\underbrace{
\left[ \begin{array}{c|c}
r_1  & r_2\\
r_1' & r_2'\\ \hline
{R_{1}}' & \mathbf{0}\\ \hline
\mathbf{0} & {R_{2}}'
\end{array} \right]
}_{R'}\underbrace{
\left[ \begin{array}{c|c}
A  & ab\\ \hline
dc & B
\end{array} \right]
}_{N}=
\underbrace{
\left[ \begin{array}{c|c}
s_1  & s_2\\
s_1' & s_2'\\ \hline
{S_{1}}' & \mathbf{0}\\ \hline
\mathbf{0} & {S_{2}}'
\end{array} \right]
}_{S'}\end{equation*}
and $[R'|S']$ is the incidence matrix associated with $N$.

\noindent
{\bf Case (c):} This case is very similar to case (b). Here we have again to find an alternative binet representation of $N_1$. 
This can be obtained if we take the representation where $e_1$ is a loop in a binet representation of $N_1$. In 
this case the incidence matrix associated with a binet representation of $N_1$ can be:
\begin{equation*}
[R_1|S_1]=
\begin{array}{c}
\begin{array} {crrrc}
& \hspace{11mm} e_3 & \hspace{0.25in} e_1 & \hspace{-1.5mm} e_2 &
\end{array}\\
\left[ \begin{array}{cr|crr}
r_1  &  1  & s_1  & 0  & 1\\
r_1' & -1  & s_1' & 2  & 1\\
{R_{1}}' & \mathbf{0} & {S_{1}}' & \mathbf{0} & \mathbf{0}
\end{array} \right]
\end{array}
\end{equation*}
and w.l.o.g. we can also assume that the incidence matrix associated with the binet matrix $N_2$ is:
\begin{equation*} 
[R_2|S_2]=
\begin{array}{c}
\begin{array} {ccccc}
\hspace{-4mm} f_3  &  & \hspace{2.8mm}  f_1 & \hspace{1.5mm} f_2 &
\end{array}\\
\left[ \begin{array}{rc|rrc}
0  & r_2 & 1  & 1  & s_2\\
2  & r_2'& 1  & -1 & s_2'\\
\mathbf{0} & {R_{2}}' & \mathbf{0} & \mathbf{0} & {S_{2}}'
\end{array} \right]
\end{array}
\end{equation*}
Using the same methodology as we did in cases (a) and (b) it can be shown that for case (c) a incidence matrix associated with 
matrix $N$, i.e. such that $R'N=S'$, is:
\begin{equation*}
[R'|S']=
\left[ \begin{array}{cc|cc}
r_1  & r_2 & s_1  & s_2\\
r_1' & r_2'& s_1' & s_2'\\ 
{R_{1}}' & \mathbf{0} & {S_{1}}'  & \mathbf{0}\\
\mathbf{0} & {R_{2}}' & \mathbf{0} & {S_{2}}'
\end{array} \right]
\end{equation*}

\noindent
{\bf Case (d):} Similarly, the incidence matrix associated with $N_2$ can be:
\begin{equation*} 
[R_2|S_2]=
\begin{array}{c}
\begin{array} {ccccc}
\hspace{-2mm} f_3  &  & \hspace{5.5mm}  f_1 & \hspace{1mm} f_2 &
\end{array}\\
\left[
\begin{array}{rc|rrc}
-1      & r_2   & -1  & 0 & s_2\\
0       & r_2'  &  1  & -1  & s_2'\\
\mathbf{0}  & {R_{2}}'   & \mathbf{0} & \mathbf{0} & {S_{2}}'
\end{array} \right]
\end{array}
\end{equation*}
We can delete the third row from matrix $[R_1|S_1]$ of (\ref{eq_st1}) in order to get a binet representation of matrix $N_1$. 
Therefore, we can assume that in 
this case the incidence matrix associated with $N_1$ can be:
\begin{equation*} 
[{R_1}|{S_1}]=
\begin{array}{c}
\begin{array} {crrrc}
& \hspace{10mm} e_3 & \hspace{0.37in} e_1 & \hspace{1mm} e_2 &
\end{array}\\\left[
\begin{array}{cr|crr}
r_1        & -1 & s_1   & 0   & -1\\
r_1'       & 1  & s_1'  & -1  & 0\\
\hat{R_1}  & \mathbf{0} &\hat{S_1}   & \mathbf{0} & \mathbf{0} 
\end{array} \right]
\end{array}
\end{equation*}
Using the same methodology as we did in all the previous cases it is easy to show that a incidence matrix associated with $N$ is:
\begin{equation} \label{eq_eighty}
[R'|S']=
\left[ \begin{array}{cc|cc}
r_1  & r_2 & s_1  & s_2\\
r_1' & r_2'& s_1' & s_2'\\ 
{R_{1}}' & \mathbf{0} & {S_{1}}'  & \mathbf{0}\\
\mathbf{0} & {R_{2}}' & \mathbf{0} & {S_{2}}'
\end{array} \right]
\end{equation}
Case (e) is directly analogous to the case (a) and (b) where $f_2$ is a link and $f_2$ is a loop, respectively. 
Case (f) is directly analogous to the case (d). For this reason we omit the proof for these cases.\\
For each of the aforementioned cases it is obvious that $[R'|S']$ is a incidence matrix of a bidirected graph, since the set of columns 
of this matrix is a combination of columns 
in $[R_1|S_1]$ and $[R_2|S_2]$. The rows/columns
of $R'$ in each case are linearly independent, something that can be proved in much the same way as we did for the $R'$ 
in Lemma \ref{lem_net-3sum-net}. Alternatively, the non-singularity of $R'$ stems also from the graphical 
explanation we give in the following section. Specifically, since there is one-to-one correspondence between the 
$R'$ and the associated bidirected graph, it can be shown that the graph induced by the edges corresponding to the columns of $R'$ form a 
negative $1$-tree in the
unique bidirected graph associated with $[R'|S']$ found in each case.
\end{proof}
\noindent
\textbf{Graphical Representation of Network $\oplus_{3}$ Binet:}\\
\noindent
An illustration regarding case (a) is depicted in Figure~\ref{fig:net3bin1}, where the triangles $(f_1,f_2,f_3)$ 
and $(e_1,e_2,e_3)$ are glued together and their edges are deleted from the unified graph. In this way, we obtain a bidirected graph 
whose associated incidence matrix is the one 
given by (\ref{eq_sixty0}).
In case (b), we convert the network representation of $N_1$ to a binet representation in which $e_2$ is a loop. As described in the 
proof of Lemma~\ref{lem_net-3sum-binet}, 
this can be done by introducing an artificial link parallel to $e_2$ and then contracting it. In this way $e_1$ becomes a negative link, 
since contraction involves switching at 
the node at which $e_1$ and $e_2$ are incident, but this does not affect the gluing of $e_1$ and $f_2$ since $f_2$ is also a negative 
link because its fundamental circuit 
uses the negative link of the basic cycle. This case is illustrated in Figure \ref{fig:net3bin2}. That figure shows the alternative 
binet representation of the matrix represented by 
the directed graph in Figure \ref{fig:net3bin1}. 
For case (c) see Figure \ref{fig:net3bin3} for an illustration. To make a similar representation for $N_1$, we can 
convert $e_1$ to a loop with a contraction. The binet graph representing $N_1$ in Figure \ref{fig:net3bin3} is an 
alternative representation to the directed graph in Figure~\ref{fig:net3bin1}. 
\begin{figure}
\begin{center}
\centering
\psfrag{sum}{\footnotesize $\oplus_{3}$}
\psfrag{=}{\footnotesize $=$}
\psfrag{e1}{\footnotesize $e_1$}
\psfrag{e2}{\footnotesize $e_2$}
\psfrag{e3}{\footnotesize $e_3$}
\psfrag{f1}{\footnotesize $f_1$}
\psfrag{f2}{\footnotesize $f_2$}
\psfrag{f3}{\footnotesize $f_3$}
\psfrag{r}{\footnotesize $r$}
\psfrag{s}{\footnotesize $s$}
\psfrag{a}{\footnotesize $a$}
\psfrag{d}{\footnotesize $d$}
\includegraphics*[scale=0.20]{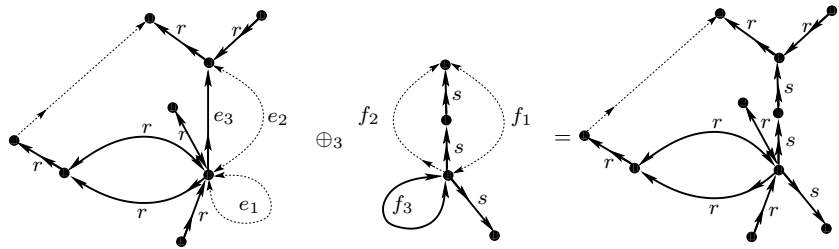}
\end{center}
\caption{The binet representation of the 3-sum of a network and a binet matrix.The case when $f_3$ is a 
loop.\label{fig:net3bin3}}
\end{figure} 
In case (d) the three edges $f_1$,$f_2$ and $f_3$ are positioned as in Figure \ref{fig:net3bin4}. We can have a 
similar position of edges $e_1,e_2,e_3$ if we delete a node that is incident to $e_1$ and $e_2$. The leftmost graph 
in Figure~\ref{fig:net3bin4} shows such a binet representation of the network matrix represented by the directed graph 
in Figure ~\ref{fig:net3bin1}. 
\begin{figure}
\begin{center}
\centering
\psfrag{sum}{\footnotesize $\oplus_{3}$}
\psfrag{=}{\footnotesize $=$}
\psfrag{e1}{\footnotesize $e_1$}
\psfrag{e2}{\footnotesize $e_2$}
\psfrag{e3}{\footnotesize $e_3$}
\psfrag{f1}{\footnotesize $f_1$}
\psfrag{f2}{\footnotesize $f_2$}
\psfrag{f3}{\footnotesize $f_3$}
\psfrag{r}{\footnotesize $r$}
\psfrag{s}{\footnotesize $s$}
\psfrag{a}{\footnotesize $a$}
\psfrag{d}{\footnotesize $d$}
\includegraphics*[scale=0.20]{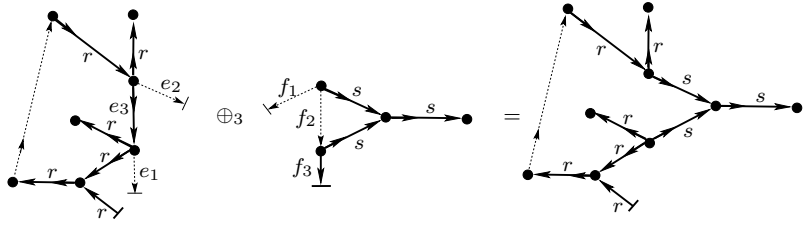}
\end{center}
\caption{The binet representation of the 3-sum of a network and a binet matrix. The case when $f_3$ is a 
half-edge.\label{fig:net3bin4}}
\end{figure}
Finally, cases (e) and (f) can be handled with the techniques described previously. If an edge among $f_1,f_2,f_3$ is a loop, 
then contract an artificial edge in the directed graph representation of $N_1$ to make the corresponding edge a loop. 
If two edges among $f_1,f_2,f_3$ are half-edges, then delete an appropriate node from the directed graph.

%%%%%%%%%%%%%%%%%%%%%%%%%%%%%%%%%%%%%%%%%%%%%%%%%%%%%%%%%%%%%%%%%%%%%%%%%%%%%%%%%%%%%%%%%%%%%%%%%%%%%%%%%%%%%%%%%%%%%%%%%%%%%%%%
\subsubsection{Binet $\oplus_3$ Network}%%%%%%%%%%%%%
A very similar analysis of the cases can be done here. The role of $e_1$, $e_2$ and $e_3$ is analogous to $f_1$, $f_2$ 
and $f_3$ as in the previous section. All the cases can be handled in much the same way, by finding a suitable alternative 
representation of $N_2$ as we did for $N_1$ in the proof of Lemma~\ref{lem_net-3sum-binet}.
\begin{lemma}
\label{lem_binet-3sum-net}
If $N_1$ is a binet matrix and $N_2$ is a network matrix such that 
\[
N_1=
\raisebox{5pt}{$
\begin{array}{r}
\vspace{3.2mm}    \\
e_3
\end{array}
\hspace{-2.5mm}
\begin{array}{c}
\begin{array}{crr}
\hspace{1mm} & e_1 & \hspace{-1.7mm} e_2 
\end{array} \\
\left[ \begin{array}{ccc}
A & a & a\\
c & 0 & 1
\end{array}
   \right],
\end{array} 
$}
N_2 = 
\raisebox{5pt}{$
\begin{array}{r}
\vspace{-5mm}
\hspace{-3.2mm}    \\
f_3
\end{array}
\hspace{-2.5mm}
\begin{array}{c}
\begin{array}{crr}
\hspace{-3.2mm} f_1 & \hspace{-1.3mm} f_2 &  
\end{array} \\
\left[ \begin{array}{ccc}
1 & 0 & b\\
d & d & B
\end{array}
   \right],
\end{array} 
$}\]
then $N= N_1 \oplus_3 N_{2}$ is a binet matrix. 
\end{lemma}

\begin{theorem}
The $k$-sum of a network (binet) matrix and a binet (network) matrix is binet ($k=1,2,3$).
\end{theorem}
\begin{proof}
The proof is similar to that of Theorem~\ref{th_net_cl} since binet matrices are also closed under duplication of columns and rows, addition of unitary rows and pivoting. 
\end{proof}

%%%%%%%%%%%%%%%%%%%%%%%%%%%%%%%%%%%%%%%%%%%%%%%%%%%%%%%%%%
\subsection{$k$-sums of Binet Matrices} \label{sec_binmat}
Here we prove that the $k$-sum ($k=2,3$) of two binet matrices is not necessarily a binet matrix. Furthermore, an analogous statement can be made for the associated 
matroids, the so-called signed-graphic matroids. Using a counterexample, we show that the $2$-sum of two binet, non-network and totally unimodular matrices, namely 
$B_1$ and $B_2$ of (\ref{eq_B1}) and (\ref{eq_B2}), is not a binet matrix. The column of $B_1$ as well as the row of $B_2$ used in our $2$-sum counterexample are indicated below: 
\begin{equation*}
B_1=
\left[ \begin{array}{c|c}
A & a
\end{array} \right]=
\left[ \begin{array}{rrrr|r}
0  & 0  & 1  & -1 & 1\\
1  & 0  & 0  &  1 & -1\\
-1 & 1  & 0  &  0 & 1\\
1  & -1 & 1  &  0 & 0\\
0  & 1  & -1 &  1 & 0
\end{array} \right],
B_2=
\left[ \begin{array}{c}
b\\ \hline
B
\end{array} \right]=
\left[ \begin{array}{rrrrr}
1  & 1  & 1  & 1 & 1\\ \hline
1  & 1  & 1  & 0 & 0\\
1 &  0  & 1  & 1 & 0\\
1  & 0 & 0  &  1 & 1\\
1  & 1  & 0 &  0 & 1
\end{array} \right]
\end{equation*}
Let $M$ be the $2$-sum of $B_1$ and $B_2$ which according to the $2$-sum definition is:
\begin{equation*}
M=
\left[ \begin{array}{c|c}
A & ab\\ \hline
\mathbf{0} & B
\end{array} \right]=
\begin{array}{c}
\begin{array}{rrrrrrrrr}
\hspace{9mm} s_1 & \hspace{1.5mm} s_2 & \hspace{1mm} s_3 & \hspace{1mm} s_4 & \hspace {1.5mm} s_5 & \hspace{1mm} 
s_6 & \hspace{1.5mm} s_7 & \hspace{1mm} s_8 & \hspace{1.5mm} s_9
\end{array}\\
\begin{array}{c}
r_1\\
r_2\\
r_3\\
r_4\\
r_5\\
r_6\\
r_7\\
r_8\\
r_9
\vspace{1mm}
\end{array}
\left[ \begin{array}{rrrr|rrrrr}
0  & 0  & 1  & -1 & 1  & 1  & 1  & 1  & 1 \\
1  & 0  & 0  &  1 & -1 & -1 & -1 & -1 & -1\\
-1 & 1  & 0  &  0 & 1  & 1  & 1  & 1  &  1\\ 
1  & -1 & 1  &  0 & 0  & 0  & 0  & 0  &  0\\
0  & 1  & -1 &  1 & 0  & 0  & 0  & 0  &  0\\ \hline
0  & 0  & 0  &  0 & 1  & 1  & 1  & 0  &  0\\
0  & 0  & 0  &  0 & 1  & 0  & 1  & 1  &  0\\
0  & 0  & 0  &  0 & 1  & 0  & 0  & 1  &  1\\
0  & 0  & 0  &  0 & 1  & 1  & 0  & 0  &  1
\end{array} \right]
\end{array}
\end{equation*} 
Assume that $M$ is a binet matrix and that $r_i$ and $s_i$ $(i=1\ldots9)$ label the basic and non-basic 
edges, respectively, in a binet representation of $M$. Matrix $M$ is integral and since it is also binet then any possible binet 
representation of $M$ up to switching 
should be one of the following two types \cite{Kotnyek:2002} (Lemmas 5.10 and 5.12):
\smallskip

\noindent {\bf Type I:} Every basic cycle is a half-edge, and all other basic edges are directed.

\noindent {\bf Type II:} There are no half-edges in the binet graph, the basis is connected and there is only one bidirected 
edge in the basis.\\
We will show that $M$ has neither of the above two representations,thereby it can not be binet. We make use of the following Lemma in~\cite{Kotnyek:2002}:
\begin{lemma} \label{lem_net_repr}
Let us suppose that a binet matrix $B$ is totally unimodular. Then it is a network matrix if and only if it has a binet 
representation in which each basic cycle is a half-edge.
\end{lemma}

\begin{lemma} \label{lem_typeI}
Matrix $M=B_1\oplus_{2}B_2$ does not have a binet representation of type I or type II.
\end{lemma}
\begin{proof}
Suppose that $M$ has a binet representation of type I. Combining the fact that $M$ is totally unimodular with 
Lemma~\ref{lem_net_repr} we have that $M$ is a network matrix. It is well-known that any submatrix of a network matrix is 
a network matrix itself (e.g. see~\cite{NemWols:1988}). $B_1$ is a submatrix of $M$ which is known to be non-network. Thus, $M$ 
can not have a binet representation of type I.

Assume that $M$ has a binet representation $\Sigma$ of type II. Let $\Sigma_R$ be the subgraph of $\Sigma$ induced by the edges 
in $R=\{r_1,\ldots,r_9\}$ ($\Sigma_R$ is also called the basis graph of $\Sigma$). Let also $C$ be the set of edges that constitute 
the unique cycle in $\Sigma_R$, i.e. $C$ is the edge set of the basic cycle of the binet graph $\Sigma$. Because of column $s_5$ of 
$M$ the subgraph of $\Sigma_R$ induced by the basic edges in $S=\{r_1,r_2,r_3,r_6,r_7,r_8,r_9\}$ is connected. Our first claim is that 
$C\subseteq{S}$. If we assume the contrary, i.e. that $C\nsubseteq{S}$, then the edges in $S$ should form a path in $\Sigma_R$. 
Moreover, observe that each non-basic edge of the set $\{s_6,s_7,s_8,s_9\}$ is using edges of $S$ in order to create the associated 
fundamental circuit in $\Sigma$. Combining this with the fact that the edges in $S$ induce a path of $\Sigma$, we have that
$\left[ \begin{array}{c}
ab\\
B
\end{array} \right]$
must be a network matrix. But this can not happen since this matrix contains $B_2$ as a submatrix which is not a network matrix and thus, 
our claim is true. Thus, $C\subseteq{S}$ and furthermore, since there is only one cycle in $\Sigma_R$, we have that $\{r_4,r_5\}\notin{C}$.

Let $D=\{r_1,r_2,r_3\}$ and $E=S-C=\{r_6,r_7,r_8,r_9\}$; our second claim is that $C\nsubseteq{D}$. If we assume the contrary, i.e. that 
$C\subseteq{D}$   then because of column $s_5$ of $M$ we have that the corresponding fundamental circuit in $\Sigma$ should be either a 
handcuff of type I or a handcuff of type II. However, it can not be a handcuff of type II since then a $\pm{2}$ would appear in $M$ 
(see Algorithm 1 in \cite{AppKot:2005}). Therefore, it is a handcuff of type I and thus the basic edges in $(D-C)\cup{E}$ induce a path 
in the basis graph. Thus, the edges in $E$ and one or more edges of $D$ are the parts of this path in the basis graph. 
Moreover, from the fundamental circuits of $\Sigma$ described by the columns of
$\left[
\begin{array}{c}
A\\
\mathbf{0}
\end{array} \right]
$ part of $M$
we have that the subgraph $\Sigma_T$ of $\Sigma_R$ induced by the set of edges in $T=\{r_1, r_2, r_3, r_4, r_5\}$ is connected. Observe 
now that the edges in $D$ appear in all the fundamental circuits of $\Sigma$ corresponding to the columns of 
$\left[ \begin{array}{c}
ab\\
B
\end{array} \right]$. Therefore, because of the structure of these fundamental circuits and since $\Sigma_T$ is connected we have that 
in $\Sigma_R$ the following conditions must be satisfied:(i) $r_6$ and $r_9$ are adjacent, (ii) $r_6$ and $r_7$ are adjacent, (iii) $r_7$ 
and $r_8$ are adjacent, and (iv) $r_8$ and $r_9$ are adjacent. We show now that this can not happen. Assume, w.l.o.g., that $r_9$ is on 
the right side of $r_6$ then because of (ii) $r_7$ should be put on the left side of $r_6$. Moreover, because of (iii) $r_8$ should be on 
the left side of $r_7$. But now condition (iv) can not be satisfied. Thus, our assumption that $C\subseteq{D}$ is not correct and this 
completes the proof of our second claim.

Since we have shown that  $\{r_4,r_5\}\notin{C}$ and that $C\nsubseteq{D}$ we have that $\Sigma_T$ is a tree in $\Sigma_R$. We show now 
that any two edges in $D$ do not share a common end-node. Note that the following procedure can be used in much the same way for any 
pair of edges in $D$. Specifically, suppose that $r_1$ and $r_2$ share an end-node and without loss of generality suppose 
that $r_2$ stands on the right side of $r_1$. Consider the fundamental circuits of $\Sigma$ determined by the columns of the
$\left[
\begin{array}{c}
A\\
\mathbf{0}
\end{array} \right]
$
part of $M$. Because of the columns $s_3$ and $s_4$ we have that $r_5$ stands on the left side of $r_1$. Moreover, 
because of the columns $s_1$ and $s_3$ we have that $r_4$ has a common end-node with $r_1$ and $r_2$. But now, we can not satisfy the 
fundamental circuit defined by $s_2$ because edge $r_1$ is in the middle of $r_4$ and $r_5$. Thus, we can conclude that any two edges 
of $D$ do not share a common end-node. However, we have that $\Sigma_T$ (which contains $r_4$ and $r_5$) is a tree and that the edges 
in $S$ (which does not contain $r_4$ and $r_5$) induce a connected subgraph in $\Sigma_R$ containing a basic cycle. This can only 
happen if $\Sigma_R$ contains at least two cycles. In other words, in order to find a binet graph satisfying the circuits described by 
the columns of $M$ we have that $\Sigma_R$ should contain at least two cycles. This is in contradiction with the fact that connected 
binet graphs contain at most one basic cycle in the basis graph. Therefore, $M$ does not have a binet representation of type II.
\end{proof}

In general we can state the following theorem.
\begin{theorem} \label{thrm_binet_ksum_open}
Totally unimodular binet matrices are not closed under $k$-sums for $k=2,3$.
\end{theorem}
\begin{proof}
For $k=2$ the Lemma \ref{lem_typeI} provides a counterexample. For $k=3$ it is enough to observe that for $c=0$ in the Definition~\ref{def_k-sums}
the $3$-sum of two matrices reduces to the $2$-sum of some submatrices obtained by the deletion of columns and rows. Since binet and TU matrices are closed under 
row and column deletions, the result follows. 
\end{proof} 

%%%%%%%%%%%%%%%%%%%%%%%%%%%%%%%%%%%%%%%%%%%%%%%%%%%%%%%%%%%%%%%%%%%%%%%%%%%%%%%%%%%%%%%%%%
\section{Tour Matrices}\label{sec_tour}

In this section a new class of matrices is introduced, that of tour matrices, in order to represent some important classes of matrices on 
bidirected graphs. In what follows, we also prove some elementary properties of tour matrices and show that they are closed under $k$-sums. 

\subsection{Definition and Elementary Properties}\label{subsec_tour_properties}
Let $[Q|S]$ be the incidence matrix of a bidirected graph $\Sigma$. We denote by $\Sigma(Q)$ and $\Sigma(S)$ the subgraphs of 
$\Sigma$ induced by the edges that correspond to the columns of $Q$ and $S$, respectively. A \emph{tour} in a bidirected graph is a walk 
in which no edge is repeated. A \emph{closed tour} is a tour in which the first and last node coincide or the first and last edge are 
half edges.
\begin{definition}
Let $\Sigma$ be a bidirected graph with $[Q|S]$ its incidence matrix. 
A $\{0,\pm1\}$ matrix $B$ with rows indexed by the columns of $Q$ and columns indexed by the columns of $S$, such that\\ 
1. $QB=S$ \\ 
2. $Q$ is full row rank\\
is called a tour matrix.
\end{definition}
The edges in $\Sigma(Q)$ are called \emph{prime} and the edges in $\Sigma(S)$ are called \emph{non-prime}. 
When in a bidirected graph representing a tour matrix $B$ the prime and non-prime edges are labeled, we call it a 
\emph{tour representation} or a \emph{tour graph} of $B$. 
\begin{lemma}
Let $B$ be an $m\times{n}$ tour matrix of a bidirected graph $\Sigma$ with incidence matrix $[Q|S]$ and $Q(b_i)$ 
be the set of edges in $Q$ indexed by the nonzero entries in the column $b_i$ of $B$ ($i=1,\ldots,n$). Then the subgraph induced 
by $Q(b_i)\cup{s_i}$ is a collection of closed tours in $\Sigma$, where $s_i$ is the $i$-th column of $S$.
\end{lemma}
\begin{proof}
Since $Qb_{i}-s_{i}=0$ for all $i\in({1,\ldots,n})$ and $q_j, s_i \in\{0,\pm1,\pm2\}^{n}, b_i\in\{0,\pm1\}^{n}$   for all $q_j\in{Q(b_{i})}$ 
we have that the degree of every vertex in the subgraph induced by $Q(b_i)\cup{s_i}$ is even, therefore its connected
components are Eulerian.  Thus, the subgraph induced by $Q(b_{i})\cup{s_i}$ is a collection of closed tours. 
\end{proof}       
In the following lemmas we provide some elementary operations which if applied to a tour matrix then the matrix 
obtained is also tour. 
\begin{lemma} \label{lem_switch}
If $\Sigma$ is a tour representation of a tour matrix $B$ then switching at a node of $\Sigma$ keeps $B$ unchanged.
\end{lemma}
\begin{proof}
Switching at a node $v$ in a bidirected graph $\Sigma$ is interpreted as multiplying by $-1$ the row of the incidence matrix 
$D=[Q|S]$ which corresponds 
to node $v$. Let $Q'$ and $S'$ be the matrices obtained multiplying by 
$-1$ the aforementioned row of $D$. Since $QB=S$, from matrix multiplication we also have that $Q'B=S'$. 
\end{proof}  

\begin{lemma} \label{lem_elem_oper}
Tour matrices are closed under the following operations:
\smallskip 

\noindent(a) Permuting rows or columns. 

\noindent(b) Multiplying a row or column by $-1$. 

\noindent(c) Duplicating a row or column.

\noindent(d) Deleting a row or column.

\end{lemma}

\begin{figure}
\begin{center}
\centering
\psfrag{f}{\footnotesize $f$}
\psfrag{u}{\footnotesize $u$}
\psfrag{s}{\footnotesize $s$}
\psfrag{f'}{\footnotesize $f'$}
\psfrag{t}{\footnotesize $t$}
\psfrag{(i)}{\footnotesize (i)}
\psfrag{(ii)}{\footnotesize (ii)}
\psfrag{(iii)}{\footnotesize (iii)}
\includegraphics*[scale=0.40]{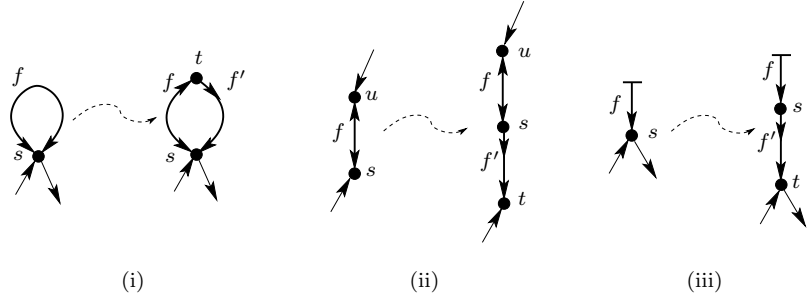}
\end{center}
\caption{The graphical equivalent of duplicating row $f$ \label{fig_duplication}} 
\end{figure}

\begin{proof}
If $B$ is a tour matrix then by definition we have that $QB=S$, where $D=[Q|S]$ is the incidence matrix of a 
bidirected graph $\Sigma$ associated with $B$. Let $B'$ be the matrix obtained by applying one of the above operations on 
$B$. We show in each case that $B'$ is a tour matrix by providing the associated incidence matrix $D'=[Q'|S']$.\\
(a) When permutation of columns takes place let $Q'=Q$ and $S'$ be the matrix obtained from $S$ by permuting the 
columns of $S$ in the same way that columns of $B$ were permuted. When permutation of rows takes place let $S'=S$ 
and $Q'$ be the matrix obtained from $Q$ by permuting its columns in the same way that rows of $B$ were permuted. 
From matrix multiplication rules we have that $Q'B'=S'$ and that $D'=[Q'|S']$ is the incidence matrix of a 
bidirected graph in both cases.  \\
(b) If row $e$ of $B$ is multiplied by $-1$ then let $Q'$ be $Q$ with column $e$ multiplied 
by $-1$ and $S'=S$. If we multiply a column $f$ of $B$ by $-1$ then let $Q'=Q$ and $S'$ be $S$ 
with column $f$ multiplied by $-1$. Obviously in both cases $B'$ is a tour matrix since from matrix multiplication 
rules we have that $Q'B'=S'$.  \\
(c) If we duplicate a column $f$ in $B$, let $Q'=Q$ and $S'$ be $S$ with column $f$ duplicated. It is easy to check
then that $B'$ satisfies the conditions of a tour matrix. 

Row duplication is a bit more involved. We have four cases corresponding to the different types of edges, and in 
each case we will alter the bidirected graph to correspond to the new tour matrix. 
If row $f$ to be duplicated is a positive loop, simply add a positive loop to any node of the signed graph. 
If the prime edge $f$ is a negative loop (see (i) in Figure~\ref{fig_duplication}), then add
a zero row $t$  in $[Q|S]$ and a zero column $f'$ in $Q$ to obtain $[Q'|S']$ and set
\[
Q'_{s,f}=Q'_{t,f}=Q_{s,f}/2 \;\;\;\mbox{and}\;\;\; Q'_{s,f'}=-Q'_{t,f'} = Q_{s,f}/2.
\]
If the prime edge $f$ is a link (see (ii) in Figure~\ref{fig_duplication}) then duplicate row $s$ in $[Q|S]$ to create
a new row $t$, and make all the elements of row $s$ zero except the element in position $f$. In row $t$ make the element 
in position $f$ zero. Finally add a new column $f'$ in $Q'$ and set
\[
Q'_{t,f'} = -Q'_{s,f'} = Q_{s,f}. 
\]
Finally, if the prime edge $f$ is a half-edge (see (iii) in Figure~\ref{fig_duplication}) then 
then add a zero row $t$  in $[Q|S]$ and a zero column $f'$ in $Q$ to obtain $[Q'|S']$ and set
\[
Q'_{s,f} = -Q'_{s,f'} = Q'_{t,f'} = Q_{s,f}.
\]
In all cases, the matrix $[Q'|S']$ is the incidence matrix of a bidirected graph by construction, and $Q'B'=S'$.  \\
(d) Deletion of a column in a tour matrix is simply the deletion of the corresponding non-prime edge in the
corresponding bidirected graph. Deletion of a row $f$, differs according the type of the corresponding prime edge $f$. 
If $f$ is a positive loop, or a link, then contract $f$ in the bidirected graph. If $f$ is a negative loop  then
make all adjacent links to the end-node of $f$ half-edges adjacent  to their other end-node, while all adjacent loops
and half edges become positive loops at some other arbitrary node, and delete $f$ and its end-node. In all cases it
is easy to verify that the new bidirected graph corresponds to the tour matrix with a column(row) deleted.  
\end{proof}
We should note here that multiplying a row (column) by -1 in a tour matrix, graphically is equivalent to reversing 
the direction of the corresponding prime (respectively non-prime) edge in the associated bidirected graph. On the 
other hand, duplicating a column amounts to creating a parallel non-prime edge
to the tour graph. 

Given a matrix 
$\left[ \begin{array}{cc}
1 & c \\
b & D 
\end{array}
   \right]$
in $\mathbb{R}$, a {\em pivot} is the matrix 
$\left[ \begin{array}{cc}
-1 & c \\
b & D-bc 
\end{array}
   \right]$
(see~\cite{Schrijver:98}). 
\begin{lemma}\label{lem_tour_pivot}
Totally unimodular tour matrices are closed under pivoting.
\end{lemma}
\begin{proof}
Let 
$T= \left[ \begin{array}{cc}
1 & c \\
b & D 
\end{array}
   \right]$
be a totally unimodular tour matrix associated with a bidirected graph $\Sigma$ with incidence matrix $[f\; Q \;|\; e\; S]$. By definition 
\begin{equation}\label{eq_20}
[f \; Q] T = [e \; S]
\end{equation}
and the columns $f$ and $e$ correspond to the prime and non-prime edges respectively. Consider the
bidirected graph $\Sigma^{'}$ with incidence matrix $[e\;Q \; | \; -f \; S]$, that is $\Sigma$ with edge $f$ having its endpoints reversed in sign. 
We will show that matrix 
$ B =  \left[ \begin{array}{cc}
-1 & c \\
b & D-bc 
\end{array}
   \right]$
is a tour matrix associated with $\Sigma^{'}$.

Initially let us show that 
\begin{equation}\label{eq_21}
[e \; Q] B = [-f \; S]
\end{equation}
We know from (\ref{eq_20}) that $f+\sum_{i} b_{i}q_{i} = e$, where $q_{i}$ is the
$i^{th}$ column of $Q$. Therefore 
\begin{equation}\label{eq_22}
-f = -e + \sum_{i} b_{i}q_{i}, 
\end{equation}
which shows that the first column of $B$ is a collection of 
tours in $\Sigma'$.  Take any other column $j$ of $B$. If $c_{j}=0$ the relationship (\ref{eq_21})  follows. If $c_{j}=+1$ then
we know from  (\ref{eq_20}) that 
\[
f + \sum_{i} d_{ij} q_{i} = s_{j}, 
\]
and the corresponding product in (\ref{eq_21}) will be 
\[
e + \sum_{i}(d_{ij}-b_{i})q_{i}.
\]
Partition the indices of the differences in the above summation into three sets: $I_{1}$ which corresponds to indices where 
both $d_{ij}, b_{i} \neq 0$, $I_{2}$ where $d_{ij}\neq 0$ and $b_{i}=0$ and $I_{3}$ where $d_{ij}=0$ and $b_{i}\neq 0$. 
Replacing $e$ by (\ref{eq_22}) we have 
\begin{eqnarray*}
e + \sum_{i}(d_{ij}-b_{i})q_{i} & = & f+\sum_{i} b_{i}q_{i} + \sum_{i\in I_{1}}(d_{ij}-b_{i})q_{i} + \sum_{i\in I_{2}}d_{ij}q_{i}
- \sum_{i\in I_{3}}b_{i}q_{i} \\
& = & f+\sum_{i} d_{ij}q_{i} = s_{j} 
\end{eqnarray*}
Similarly for the case where $c_{j}=-1$ (or alternatively use (b) of Lemma~\ref{lem_elem_oper}). 

Given that totally unimodular matrices  are closed under pivoting,  $B$ will be a $\{0,\pm1\}$ matrix. 
\end{proof}

\begin{lemma} \label{lem_nettour}
Network matrices are tour matrices.
\end{lemma}
\begin{proof}
Consider a network matrix $N\in\{0,\pm1\}$ of a directed graph $G$ with incidence matrix $[R|S]$. We will show that $N$ can be viewed as a 
binet matrix by providing a binet representation of it.

Let $e$ be any column of $S$.  In what follows we will show that there exists a binet representation in which 
edge $e$ is a loop at any one of its endpoints. View the graph $G$ as a bidirected graph $\Sigma$ with only 
positive links. Add a negative link $f$ parallel to $e$ and as a result we 
have that the binet matrix associated with $\Sigma$ is equal to the original network matrix $N$ plus an all-zero 
row. Deleting this all-zero row  we get the original matrix $N$, while the  equivalent graphical operation would be the
contraction of edge $f$. Contraction of $f$ involves switching at one end-node of $f$ (say at $v$), since $f$ is a 
negative link. This way $e$ becomes a negative loop (see Figure \ref{fig:netcycle}).
\begin{figure}[h]
\centering
\psfrag{e}{\footnotesize $e$}
\psfrag{f}{\footnotesize $f$}
\psfrag{v}{\footnotesize $v$}
\includegraphics*[scale=.15]{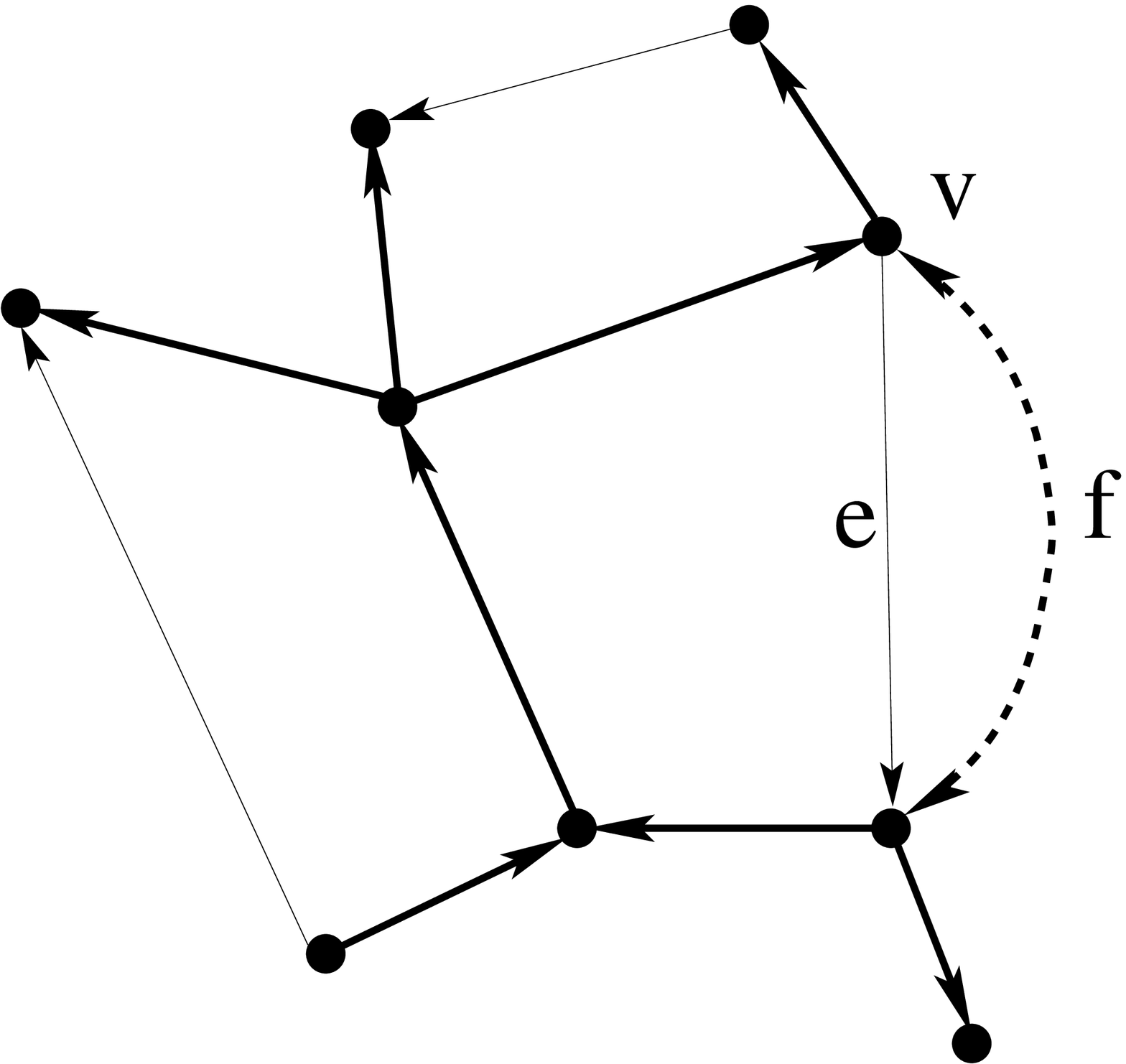} \hspace{0.2in}
\includegraphics*[scale=.15]{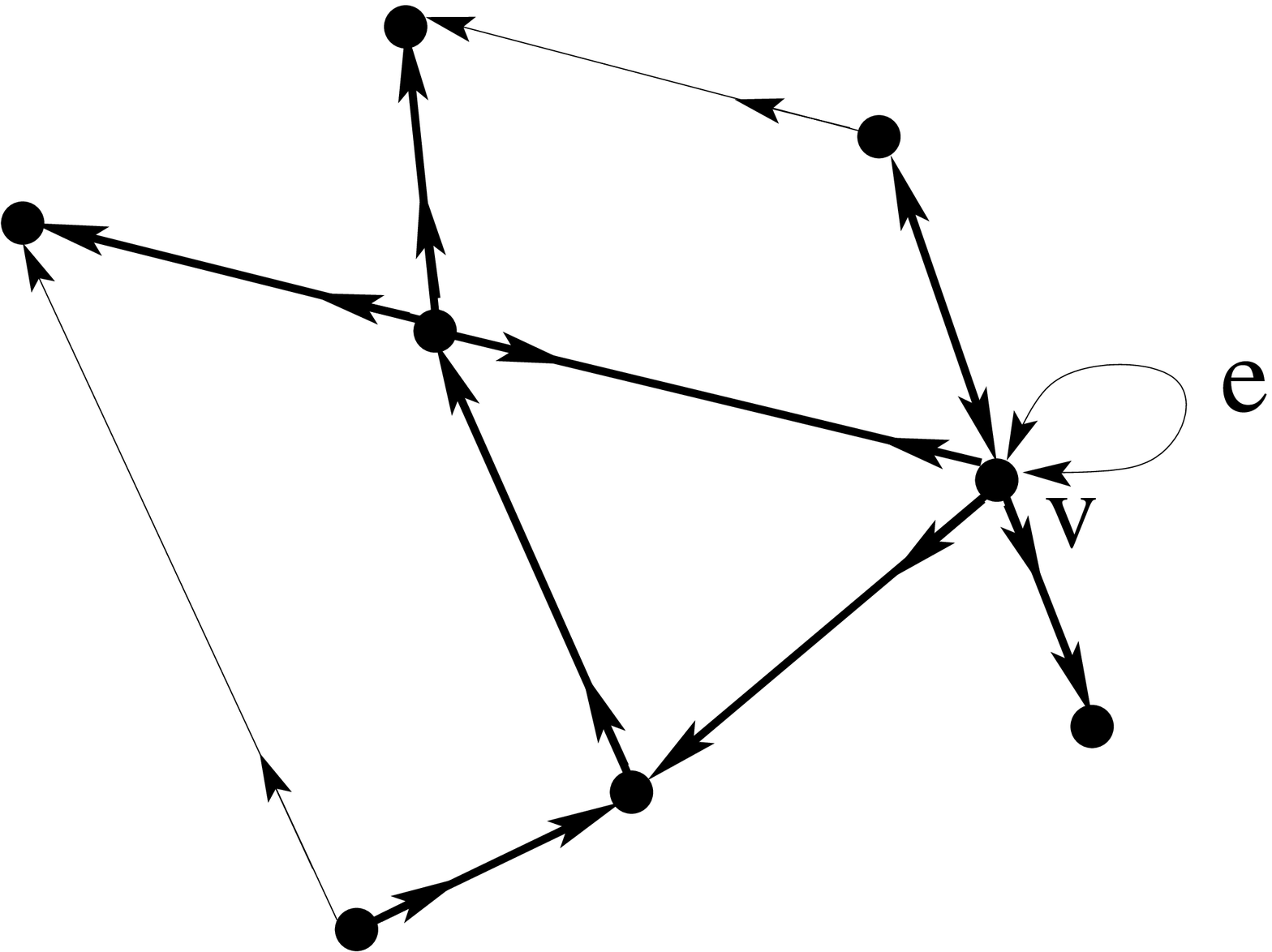} %\hspace*{\fill}
\caption{Inserting a negative edge $f$, and then contracting it by switching at $v$.\label{fig:netcycle}}
\end{figure}

In matrix terms, we have that starting from 
$[R|S]=\begin{array}{r}
v\\
\\ 
\end{array}\!\!\!\!\!
\begin{array}{c}
\hspace{3mm} e\\
\left[ \begin{array}{c|rc}
 & -1 & \\
R    &  1 & \bar{S}\\
 & \mathbf{0} &  
\end{array} \right]
\end{array}
$ 
by the aforementioned procedure we obtain $[R'|S']=
\begin{array}{r}
v\\
\\ 
\end{array}\!\!\!\!\!
\begin{array}{c}
\hspace{1mm} e\\
\left[ \begin{array}{c|cr}
    & 2 & \\
\hat{R} & 0 & \hat{S} \\
    &\mathbf{0} & 
\end{array} \right]
\end{array}$,
where $R'N=S'$ and $[R'|S']$ is a incidence matrix associated with a binet representation of $N$. 
Therefore we have found a bidirected graph $[R'|S']$ where $R'N=S'$, and $R'$ is full-row
rank. 
\end{proof}
Furthermore, it is known that any binet matrix which is TU and non-network should have a binet representation $\Sigma$ that does not 
contain half-edges (see Lemma~22 and Theorem~24 in~\cite{AppKot:2005}). Therefore we can state the 
following corollary. 
\begin{corollary} \label{cor_tubin}
TU binet matrices are tour matrices.
\end{corollary}
From Corollary~\ref{cor_tubin}, it is evident that $B_1$ and $B_2$ are tour matrices. Combining 
this with Lemma~\ref{lem_nettour}, Theorem~\ref{Seymour_matrix} and the fact that zero columns are preserved,  
we have the following.
\begin{corollary} \label{col_tutour}
All the building blocks of TU matrices  are tour matrices and their transposes.
\end{corollary}   

%%%%%%%%%%%%%%%%%%%%%%%%%%%%%%%%%%%%%%%%%%%%%%%%%%%%%%%%%%%%%%%%%%%%%%%%%%%%%%%%%%%%%%%

\subsection{Bidirected Graph Representation of TU matrices}
In this section we will show that all TU matrices have a bidirected graph representation since they are a subclass
of tour matrices. This is illustrated in the following ``pathological'' case by the usage of positive loops, 
which in general allow a somewhat arbitrary insertion of prime edges and thereby rows in a given matrix.  
\begin{theorem}\label{thrm_tu=tour}
All TU matrices are tour matrices.
\end{theorem}
\begin{proof}
Let $B\in \{0,\pm1\}^n\times m$ be a totally unimodular matrix.
By Ghouila-Houri characterisation of TU matrices (see~\cite{Ghouila:62}), we have that 
there exists a vector $x^{T}\in \{\pm1\}^{n}$ such that $x^{T}B = y^{T} \in \{0,\pm1\}^{n}$; that is multiplying the rows
by $\pm1$ the resulting matrix has columns which sum up to $\{0,\pm1\}$. Therefore we can have
$
\left[ \begin{array}{c}
x^{T} \\
x^{T}
\end{array}
   \right]
B = \left[ \begin{array}{c}
y^{T} \\
y^{T}
\end{array}
   \right]
$
and 
$
[Q | S] = \left[ \begin{array}{cc}
x^{T} & y^{T} \\
x^{T} & y^{T}
\end{array}
   \right],
$
is the incidence matrix of a bidirected graph since the sum of each column is less or equal to $|2|$. If the first column of
$Q$ is $\left[ \begin{array}{c} -1 \\ -1 \end{array} \right]$ replace it with $\left[ \begin{array}{c} -2 \\ 0 \end{array} \right]$, 
while if it is $\left[ \begin{array}{c} 1 \\ 1 \end{array} \right]$ replace it with $\left[ \begin{array}{c} 2 \\ 0 \end{array} \right]$ 
to obtain a new matrix $Q'$, and set $S'=Q'B$. Then $[Q'|S']$ is also the incidence matrix of a bidirected graph with $B$ its tour matrix. 
\end{proof}
However a tour matrix may have multiple bidirected graph representations, and in 
the proof of Theorem~\ref{thrm_tu=tour} the bidirected graph so constructed does not have enough 
structural information with respect to the linear independence of the columns of the associated matrix. We know from 
Seymour's decomposition Theorem~\ref{Seymour_matrix} that a TU matrix is composed by $k$-sums from matrices which do
have a bidirected graph representation, therefore in view of Corollary~\ref{col_tutour} there must exist a richer in structure 
bidirected graph representation. Moreover, the building blocks in the $k$-sum composition do have bidirected graphs which 
do not have positive loops. 
In order to obtain this representation, we have to examine the way the $k$-sum operations 
behave on tour matrices.

\subsubsection{The $k$-sum Operations on Tour Matrices}\label{subsubsec_ksums_tour}
In what follows we present results on the $k$-sums of tour matrices. The case of only $3$-sum will be shown as we did in the previous 
sections, since the other sum operations could be reduced to it by the addition of unitary rows and duplication of columns.

\begin{lemma} \label{lem_3-sum_tour_operations}
If $K, L$ are tour matrices, then there exist tour matrices $K', L'$ such that $K \oplus_3 L$ is a row submatrix
of $K' \oplus_3 L'$ where the connecting elements are all positive links.
\end{lemma}
\begin{proof}
Let 
\[
K=
\raisebox{5pt}{$
\begin{array}{r}
\vspace{3.2mm}    \\
e_3
\end{array}
\hspace{-2.5mm}
\begin{array}{c}
\begin{array}{crr}
\hspace{1mm} & e_1 & \hspace{-1.7mm} e_2 
\end{array} \\
\left[ \begin{array}{ccc}
A & a & a\\
c & 0 & 1
\end{array}
   \right],
\end{array} 
$}
L = 
\raisebox{5pt}{$
\begin{array}{r}
\vspace{-5mm}
\hspace{-3.2mm}    \\
f_3
\end{array}
\hspace{-2.5mm}
\begin{array}{c}
\begin{array}{crr}
\hspace{-3.2mm} f_1 & \hspace{-1.3mm} f_2 &  
\end{array} \\
\left[ \begin{array}{ccc}
1 & 0 & b\\
d & d & B
\end{array}
   \right].
\end{array} 
$}\]
For all possible edge type configurations of $f_1,f_2$ 
and $f_3$ we will apply graphical operations on the tour graph of $L$, so that the resulting graph will be 
the tour graph of a tour matrix $L'$ that will contain $L$ as a submatrix. 

\noindent
{\bf Case (a):}
Consider the case where $f_1$ is a negative loop , $f_2$ is a negative link  and $f_3$ is a positive link. Because of the first two columns 
of $L$ we have that these edges must be of the following form: $f_1=\{v,v\}$, $f_2=\{u,v\}$ and $f_3=\{u,v\}$ (see Figure~\ref{fig_case_b}). 
The graphical operation is the following: we split the end-node $v$ of $f_1$ into two nodes $v_1$ and $v_2$ and add a new  basic 
positive link $f'=\{v_1,v_2\}$ . In the new bidirected graph  $f_1, f_2$ are negative links, and $f_3$ is a positive link, while for all 
other edges having end-node $v$ we replace $v$ by $v_1$. 
Up to switchings, the tour matrix $L'$ associated with this  graph is:
\raisebox{5pt}{$
\begin{array}{r}
\vspace{-1mm}
\hspace{-3.2mm}    \\
f_3\\
\\
f'
\end{array}
\hspace{-2.5mm}
\begin{array}{c}
\begin{array}{crr}
\hspace{-3.2mm} f_1 & \hspace{-1.3mm} f_2 &  
\end{array} \\
\left[ \begin{array}{ccc}
1 & 0 & b\\
d & d & B \\
1 & 1 & \mathbf{0}
\end{array}
   \right],
\end{array} 
$}
where the connecting elements $f_1,f_2$ and $f_3$ are all positive links.
\begin{figure}[h]
\begin{center}
\centering
\psfrag{=}{\footnotesize $\rightarrow$}
\psfrag{f1}{\footnotesize $f_1$}
\psfrag{f2}{\footnotesize $f_2$}
\psfrag{f3}{\footnotesize $f_3$}
\psfrag{f'}{\footnotesize $f'$}
\psfrag{v}{\footnotesize $v$}
\psfrag{u}{\footnotesize $u$}
\psfrag{v1}{\footnotesize $v_1$}
\psfrag{v2}{\footnotesize $v_2$}
\includegraphics*[scale=0.25]{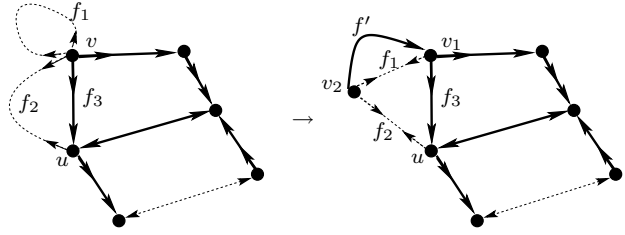}
\end{center}
\caption{$f_1$ negative loop , $f_2$ negative link  and $f_3$ positive link.\label{fig_case_b}}
\end{figure}

\noindent
{\bf Case (b):}
Let now  $f_1=\{u\}, f_3=\{v\}$ be half-edges and $f_2=\{u,v\}$ a positive link (see Figure~\ref{fig_case_d}). 
The graphical operation in this case is the following: add a new vertex $w$ in the bidirected graph, 
replace the half-edges  $f_1$ and $f_3$ by positive links  $f_1=\{u,w\}$ and $f_3=\{v,w\}$, and add a negative loop $f'=\{w,w\}$. 
The new tour matrix $L'$ associated with this graph will be:
\raisebox{5pt}{$
\begin{array}{r}
\vspace{-1mm}
\hspace{-3.2mm}    \\
f_3\\
\\
f'
\end{array}
\hspace{-2.5mm}
\begin{array}{c}
\begin{array}{crr}
\hspace{-3.2mm} f_1 & \hspace{-1.3mm} f_2 &  
\end{array} \\
\left[ \begin{array}{ccc}
1 & 0 & b\\
d & d & B \\
0 & 0 & b
\end{array}
   \right].
\end{array} 
$}
\begin{figure}[h]
\begin{center}
\centering
\psfrag{=}{\footnotesize $\rightarrow$}
\psfrag{f1}{\footnotesize $f_1$}
\psfrag{f2}{\footnotesize $f_2$}
\psfrag{f3}{\footnotesize $f_3$}
\psfrag{f'}{\footnotesize $f'$}
\psfrag{v}{\footnotesize $v$}
\psfrag{u}{\footnotesize $u$}
\psfrag{v1}{\footnotesize $v_1$}
\psfrag{v2}{\footnotesize $v_2$}
\psfrag{w}{\footnotesize $w$}
\includegraphics*[scale=0.25]{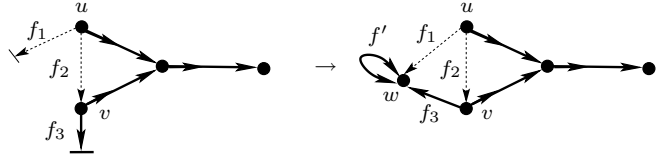}
\end{center}
\caption{$f_1, f_3$ half-edges and $f_2$ positive link.\label{fig_case_d}}
\end{figure}

\noindent
{\bf Case (c):}
For the case where $f_1,f_2$ are negative loops and $f_3$ a positive loop the graphical operation is similar to the ones
described previously, and is depicted in  Figure~\ref{fig_case_e}. 
\begin{figure}[h]
\begin{center}
\centering
\psfrag{=}{\footnotesize $\rightarrow$}
\psfrag{f1}{\footnotesize $f_1$}
\psfrag{f2}{\footnotesize $f_2$}
\psfrag{f3}{\footnotesize $f_3$}
\psfrag{f'}{\footnotesize $f'$}
\psfrag{f''}{\footnotesize $f''$}
\includegraphics*[scale=0.25]{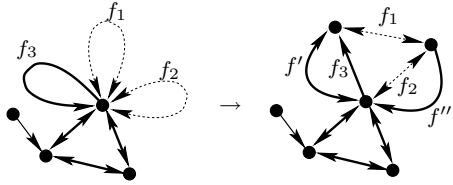}
\end{center}
\caption{$f_1, f_2$ negative loops and $f_3$ positive loop.\label{fig_case_e}}
\end{figure}
The new tour matrix $L'$ associated with the so constructed graph will be:
\raisebox{5pt}{$
\begin{array}{r}
\vspace{-1mm}
\hspace{-3.2mm}    \\
f_3\\
\\
f'  \\
f''
\end{array}
\hspace{-2.5mm}
\begin{array}{c}
\begin{array}{crr}
\hspace{-3.2mm} f_1 & \hspace{-1.3mm} f_2 &  
\end{array} \\
\left[ \begin{array}{ccc}
1 & 0 & b\\
d & d & B \\
0 & 0 & b \\
1 & 1 & \mathbf{0}
\end{array}
   \right].
\end{array} 
$}

\noindent
{\bf Case (d):}
The case where $f_1$ is a negative link, $f_2$ a negative loop and $f_3$ a half-edge, can be easily verified that is not
possible, due to the structure of $L$.

It is straightforward to show that all possible edge type configurations for the connecting edges of $L$, 
fall into one of the above described cases where the new tour matrix $L'$ will contain either a row $f'$ or $f''$ or both. 
Applying the above graphical operations and switchings on both $K$ and $L$, we can therefore obtain $K'$ and $L'$ were the
connecting elements of $K' \oplus_3 L'$ are positive links $e_{1},e_2, e_3$ and $f_1,f_2,f_3$, while the matrix $K' \oplus_3 L'$
contains $K \oplus_3 L$ as a row submatrix. 
\end{proof}

\begin{lemma} \label{lem_3-sum_tour}
If $K, L$ are tour matrices such that 
\[
K=
\raisebox{5pt}{$
\begin{array}{r}
\vspace{3.2mm}    \\
e_3
\end{array}
\hspace{-2.5mm}
\begin{array}{c}
\begin{array}{crr}
\hspace{1mm} & e_1 & \hspace{-1.7mm} e_2 
\end{array} \\
\left[ \begin{array}{ccc}
A & a & a\\
c & 0 & 1
\end{array}
   \right],
\end{array} 
$}
L = 
\raisebox{5pt}{$
\begin{array}{r}
\vspace{-5mm}
\hspace{-3.2mm}    \\
f_3
\end{array}
\hspace{-2.5mm}
\begin{array}{c}
\begin{array}{crr}
\hspace{-3.2mm} f_1 & \hspace{-1.3mm} f_2 &  
\end{array} \\
\left[ \begin{array}{ccc}
1 & 0 & b\\
d & d & B
\end{array}
   \right],
\end{array} 
$}\]
then $M= K \oplus_3 L$ is a tour matrix. 
\end{lemma}

\begin{proof}
Let $D_1=[Q_1|S_1]$ and $D_2=[Q_2|S_2]$ be incidence matrices associated with $K$ and $L$  and  
$\Sigma(D_1)$ and $\Sigma(D_2)$ be the associated tour graphs.
By Lemma~\ref{lem_3-sum_tour_operations} and (d) of Lemma~\ref{lem_elem_oper}, we can assume that the connecting elements 
$e_{1},e_2, e_3$ and $f_1,f_2,f_3$ are all positive links in the tour graphs. 

By Lemma \ref{lem_elem_oper}, the incidence matrices $D_1$ and $D_2$ associated with $K$ and $L$
can have the following form:
\begin{equation} \label{eq_hundred} 
[Q_1|S_1]=
\begin{array}{c}
\begin{array} {crrrc}
& \hspace{0.37in} e_3 & \hspace{0.34in} e_1 & \hspace{1mm} e_2 &
\end{array}\\
\left[
\begin{array}{cr|crr}
q_1        & -1 & s_1   & 0   & -1\\
q_1'       & 1  & s_1'  & -1  & 0\\
q_1''      & 0  & s_1'' & 1   & 1\\
{Q_{1}}'  & \mathbf{0} &{S_{1}}'   & \mathbf{0} & \mathbf{0} 
\end{array} \right]
\end{array}
\hspace{-4mm}
\begin{array}{l}
u  \\
v  \\
y  \\
\vspace{-4mm}
\end{array}, \quad 
[Q_2|S_2]=
\begin{array}{r}
u'  \\
v'  \\
y'  \\

\end{array}
\hspace{-3mm}
\begin{array}{c}
\begin{array} {ccccc}
\hspace{-2mm} f_3  &  & \hspace{5.5mm}  f_1 & \hspace{1mm} f_2 &
\end{array}\\
\left[
\begin{array}{rc|rrc}
0      & q_2   & -1  & -1 & s_2\\
-1     & q_2'  &  0  & 1  & s_2'\\
1      & q_2'' &  1  & 0  & s_2''\\
\mathbf{0}  & {Q_{2}}'   & \mathbf{0} & \mathbf{0} & {S_{2}}'
\end{array} \right]
\end{array}
\end{equation}
where $\mathbf{0}$ is a vector or matrix of zeroes of appropriate size, $q_{i},q_{i}',q_{i}'',s_{i},s{i}'$ and $s_{i}''$ are row
vectors and ${Q_{i}}', {S_{i}}'$ are matrices of appropriate size $(i=1,2)$. Also, $u,v$ and $y$
label the three first rows of $D_1$ and consequently the corresponding nodes of $\Sigma(D_1)$. Similarly, $u',v'$, $y'$ label the 
first three rows of $D_2$ and the corresponding nodes of $\Sigma(D_2)$. We have that the following equations hold:
\begin{equation} \label{eq_hundredtwo}
Q_1K=S_1, \qquad
Q_2L=S_2
\end{equation}     
For $K$ using (\ref{eq_hundred}) and (\ref{eq_hundredtwo}) we have that:
\begin{equation*} 
\left[
\begin{array}{c|r}
q_1  &  -1\\ 
q_1' &   1\\
q_1''&   0\\ \hline
{Q_{1}}' & \mathbf{0}
\end{array} \right]
\left[
\begin{array}{c|c|c}
A  &  a  & a \\ \hline
c  &  0  & 1
\end{array}\right]=
\left[ \begin{array}{c|r|r}
s_1       &   0  & -1\\
s_1'      &  -1  &  0\\
s_1''     &   1  &  1\\ \hline
{S_{1}}' & \mathbf{0} & \mathbf{0}
\end{array} \right]
\end{equation*}
From the above equation we take the following equations:
\begin{gather} \label{eq_hundredfive}
\left[ \begin{array}{c}
q_1\\
q_1'\\
q_1''
\end{array} \right]A+
\left[ \begin{array}{r}
-1\\
1\\
0
\end{array} \right]c=
\left[ \begin{array}{c}
s_1\\
s_1'\\
s_1''
\end{array} \right], \quad
\left[ \begin{array}{c}
q_1\\
q_1'\\
q_1''
\end{array} \right]a =
\left[ \begin{array}{r}
0\\
-1\\
1
\end{array} \right], \quad \nonumber \\
\left[ \begin{array}{c}
q_1\\
q_1'\\
q_1''\\
\end{array} \right]a+
\left[\begin{array}{r}
-1\\
1\\
0
\end{array} \right]=
\left[ \begin{array}{r}
-1\\
0\\
1
\end{array} \right], \quad
{Q_{1}}'A={S_{1}}', \quad
{Q_{1}}'a=\mathbf{0}
\end{gather}
Similarly, for $L$ using (\ref{eq_hundred}) and (\ref{eq_hundredtwo}) we have that:
\begin{equation*}
\left[ \begin{array}{r|c}
0    &   q_2\\
-1   &   q_2'\\
1    &   q_2''\\ \hline
\mathbf{0} & {Q_{2}}'
\end{array} \right]
\left[ \begin{array}{c|c|c}
1    &  0   &  b\\ \hline
d    &  d   &  B
\end{array} \right]=
\left[ \begin{array}{r|r|c}
-1   &  -1  & s_2\\
0    &   1  & s_2'\\
1    &   0  & s_2''\\ \hline
\mathbf{0} & \mathbf{0} & {S_{2}}'
\end{array} \right]
\end{equation*}
From the above equation we take the following equations:
\begin{gather} \label{eq_hundredeleven}
\left[ \begin{array}{r}
0\\
-1\\
1
\end{array} \right]+
\left[ \begin{array}{c}
q_2\\
q_2'\\
q_2''
\end{array} \right]d=
\left[ \begin{array}{r}
-1\\
0\\
1
\end{array} \right], \quad
\left[ \begin{array}{c}
q_2\\
q_2'\\
q_2''
\end{array} \right]d=
\left[ \begin{array}{r}
-1\\
1\\
0
\end{array} \right], \quad \nonumber \\
\left[ \begin{array}{r}
0\\
-1\\
1
\end{array} \right]b+
\left[ \begin{array}{c}
q_2\\
q_2'\\
q_2''
\end{array} \right]B=
\left[ \begin{array}{c}
s_2\\
s_2'\\
s_2''
\end{array} \right], \quad
{Q_{2}}'d=\mathbf{0}, \quad
{Q_{2}}'B={S_{2}}' \quad
\end{gather}
Using block matrix multiplication and equations in (23) and (24), it is easy to show that the following equation holds:
\begin{equation} \label{eq_fiftynine1}
\underbrace{
\left[ \begin{array}{c|c}
q_1  & q_2\\
q_1' & q_2'\\
q_1''& q_2''\\ \hline
{Q_{1}}' & \mathbf{0}\\ \hline
\mathbf{0} & {Q_{2}}'
\end{array} \right]
}_{Q'}\underbrace{
\left[ \begin{array}{c|c}
A  & ab\\ \hline
dc & B
\end{array} \right]
}_{M}=
\underbrace{
\left[ \begin{array}{c|c}
s_1  & s_2\\
s_1' & s_2'\\
s_2''& s_2''\\ \hline
{S_{1}}' & \mathbf{0}\\ \hline
\mathbf{0} & {S_{2}}'
\end{array} \right]
}_{S'}
\end{equation}
Clearly $D'=[Q'|S']$ is incidence matrix of a bidirected graph. 
\end{proof}
Let us examine the structure of the bidirected graph $\Sigma(D')$ so obtained, from the $k$-sum operation on tour matrices. 
From (\ref{eq_fiftynine1}) we have that $\Sigma(D')$ 
is obtained by gluing $\Sigma(D_1)$ and $\Sigma(D_2)$ such that $u$ and $u'$, $v$ and $v'$, $y$ and $y'$ become single 
nodes $u$, $v$ and $y$, respectively, and deleting edges $e_1$, $e_2$, $e_3$, $f_1$, $f_2$ and $f_3$ from the 
unified graph. In other words, this can also be seen as gluing together the $\Sigma(D_1)$ and $\Sigma(D_2)$ along the 
triangles $(e_1, e_2, e_3)$ and $(f_1, f_2, f_3)$ so that $e_1$ meets $f_3$, $e_2$ meets $f_1$ and $e_3$ meets $f_2$ 
and deleting the glued triangle from the unified graph. Obviously, we can say that in $\Sigma(D')$ the edge $e_3$ which 
was deleted is substituted by the tour associated with $f_2$ in $\Sigma(D_2)$ and that the $f_3$ which was deleted is substituted 
by the tour associated with $e_1$ in $\Sigma(D_1)$. 
Therefore, now any tour that used $e_3$ will instead go through the tour associated with $f_2$ giving rise to the 
non-zero part of $dc$ in $K\oplus_{3}L$. The tours that went through $f_3$ use the tour of $e_1$ in the unified graph, 
as determined by the $ab$ part of $K\oplus_{3}L$. All other tours remain unchanged, as it is expressed by the fact 
that if $c$ or $b$ had a zero element then $dc$ or $ab$ has an all-zero column in the same position.

From Lemmata~\ref{lem_3-sum_tour} and~\ref{lem_tour_pivot} and the fact that $1$-, and $2$-sum operations are special cases of the 
$3$-sum operation  we obtain the following theorem:
\begin{theorem} \label{th_tourksums}
Totally unimodular tour matrices are closed under $k$-sums for $k=1,2,3$.
\end{theorem}

\subsubsection{Graph Algorithm}\label{subsubsec_algorithm}
We are now ready to present an algorithm which given a totally unimodular matrix $N$ will construct a bidirected graph
$\Sigma$ or equivalently an incidence matrix, where the columns in $N$  represent collection of closed tours.

\begin{itemize}

\item[1.] Given a TU matrix $N$, by Seymour's Theorem~\ref{Seymour_matrix} we can decompose it via $k$-sums
into matrices $N_1, \ldots, N_{n}$. A separation algorithm for this can be found in the book by Truemper~\cite{Truemper:98}.

\item[2.] For each matrix $N_{i}$ one of the following cases will be true:

\begin{itemize}
\item[2.1] Check whether $N_{i}$ is a network matrix, and if so construct the associated  incidence matrix $D_{\Sigma_{i}}$. This
can be done by the Tutte's recognition algorithm which results from his  decomposition theory for graphic 
matroids~\cite{Tutte:1960,BixCunn:1980}. 

\item[2.2] Check whether $N_{i}$ is a binet matrix, and if so construct the associated incidence matrix $D_{\Sigma_{i}}$. Similarly
with step 2.1, a decomposition theory for binary signed graphic matroids given in~\cite{AppPapPit:2008}  can be used
in this step. Alternatively one can also use the algorithm given in~\cite{Musitelli:2007}.

\item[2.3] If neither of the above cases is true, then $N_{i}$ is the transpose of a network matrix which is not
binet. In this case construct the bidirected graph representation given in the proof of Theorem \ref{thrm_tu=tour}.
\end{itemize}

\item[3] Starting from the incidence matrices $D_{\Sigma_{i}}$ resulted from step 2 and the $k$-sum decomposition indicated in step 1, 
compose the incidence matrix of $N$ using the matrix operations so defined in the constructive proofs of Lemmata~\ref{lem_net-3sum-net}, 
\ref{lem_net-3sum-binet} and \ref{lem_3-sum_tour}.

\end{itemize}

All of the above steps can be performed in polynomial time with respect to the size of the matrix $N$.

The fact that case 2.3 in the above algorithm is possible, that is the existence of a transpose of a network matrix which is
not binet, is verified by a  recent work of Slilaty~\cite{Slilaty:2005b} where he identifies a set of 29 excluded minors for a 
cographic matroid to be signed graphic. Examination of the aforementioned excluded minors, reveals that all are a 2- or 
3-sum of two binet matrices without positive loops, therefore by Lemma~\ref{lem_3-sum_tour}, tour matrices with a bidirected
graph representation without positive loops. However, we were unable to generalise this to an arbitrary non-binet 
transpose of a network matrix, therefore we use the trivial bidirected graph representation given in the proof of Theorem \ref{thrm_tu=tour}.

\section{Concluding Remarks}
Totally unimodular matrices characterise a class of well solved integer programming problems, due to the integrality
property of the associated polyhedron. In this paper we exploit the decomposition theorem of Seymour for totally unimodular matrices,
and provide a graphical representation for every such matrix in a bidirected graph, such that the structural information
of the decomposition building blocks is mostly retained.  In order to do this, we
examine the effect of the $k$-sum operations on network matrices, their transposes and binet matrices, and show that
the aforementioned classes of matrices are not closed under these composition operations. A new, more general, 
class of matrices is introduced called tour matrices, which is proved to be closed under $k$-sums, and it has an
associated bidirected graph representation in the sense that the columns of a tour matrix represent a collection of closed tours.

\section{Acknowledgements}
The authors wish to thank Thomas Zaslavsky  and the two anonymous referees whose comments helped to improve 
and enhance the  presentation of the results in this paper.

\bibliographystyle{plain}

\end{document}